\documentclass[a4paper,11pt]{article}

\usepackage[T1]{fontenc}
\usepackage{lmodern}

\usepackage[latin1]{inputenc}
\usepackage{amsmath}
\usepackage{amsthm}
\usepackage{amssymb}
\usepackage{mathrsfs}
\usepackage{graphicx}
\usepackage[all]{xy}
\usepackage{hyperref}

\usepackage{abstract} 

\newtheorem{thm}{Theorem}[section]
\newtheorem{cor}[thm]{Corollary}

\newtheorem{fact}[thm]{Fact}
\newtheorem{lemma}[thm]{Lemma}
\newtheorem{prop}[thm]{Proposition}

\theoremstyle{definition}
\newtheorem{definition}[thm]{Definition}
\newtheorem{ex}[thm]{Example}
\newtheorem{remark}[thm]{Remark}

\title{Coning-off CAT(0) cube complexes}
\date{\today}
\author{Anthony Genevois}

\begin{document}

\maketitle

\begin{abstract}
In this paper, we study the geometry of cone-offs of CAT(0) cube complexes over a family of combinatorially convex subcomplexes, with an emphasis on their Gromov-hyperbolicity. A first application gives a direct cubical proof of the characterization of the (strong) relative hyperbolicity of right-angled Coxeter groups, which is a particular case of a result due to Behrstock, Caprace and Hagen. A second application gives the acylindrical hyperbolicity of $C'(1/4)-T(4)$ small cancellation quotients of free products. 
\end{abstract}

\tableofcontents

\newpage

\section{Introduction}

\hspace{0.5cm} A fruitful method to study a given group is to make it act on some space which is ``negatively-curved'' in some sense. Following Gromov, a possibility is to define \emph{$\delta$-hyperbolic} geodesic spaces by requiring that, in any geodesic triangle, any side is included into the $\delta$-neighborhood of the union of the two other sides. Often it happens that a group $G$ admits a natural action on some geodesic space $X$, which is generally not hyperbolic however. A well-known method to produce an action on a hyperbolic space is to \emph{cone-off} this space: loosely speaking, we construct a new geodesic space $Y$ from $X$ by gluing cones over subspaces of $X$, in order to ``kill'' the non-hyperbolic subspaces of $X$ and to make $Y$ hyperbolic; thus, if this collection of subspaces is $G$-equivariant, the action $G \curvearrowright X$ naturally induces a new action $G \curvearrowright Y$. In this paper, for convenience we will use two different definitions of a cone-off; nevertheless, the spaces we obtain are quasi-isometric:

\begin{definition}\label{cone-off}
Let $X$ be a CW complex and $\mathcal{Q}$ a collection of subcomplexes. The \emph{cone-off of $X$ over $\mathcal{Q}$} is the graph obtained from $X^{(1)}$ by adding an edge between two vertices whenever they both belong to a common subcomplex of $\mathcal{Q}$. 
\end{definition}

\begin{definition}
Let $X$ be a CAT(0) cube complex and $\mathcal{Q}$ a collection of subcomplexes. The \emph{usual cone-off of $X$ over $\mathcal{Q}$} is the graph obtained from $X^{(1)}$ by adding a vertex for each subcomplex $Q \in \mathcal{Q}$ and linking it by an edge $Q$ to each vertex belonging to $Q$.
\end{definition}

In this article, we focus on the class of \emph{CAT(0) cube complexes}: roughly speaking, a CAT(0) cube complex is a geodesic space constructed by gluing cubes together, in such a way that the geodesic triangles turn out to be thiner than the Euclidean triangles. Essentially, we determine precisely when a CAT(0) cube complex is hyperbolic, in order to identify the possible obstructions to hyperbolicity, and we prove that the cone-off of a CAT(0) cube complex, killing these obstructions, turns out to be hyperbolic. More precisely, we first prove:

\begin{thm}\label{hyperbolic}
Let $X$ be a CAT(0) cube complex. The following are equivalent:
\begin{itemize}
	\item[(i)] $X$ is hyperbolic;
	\item[(ii)] the flat rectangles in $X$ are uniformly thin;
	\item[(iii)] $X$ is finite-dimensional and the grid of hyperplanes in $X$ are uniformly thin.
\end{itemize}
\end{thm}

A \emph{flat rectangle} is a combinatorially geodesic subcomplex isomorphic to some square complex $[0,a] \times [0,b]$; it is \emph{$L$-thin} if $\min(a,b) \leq L$ and \emph{$L$-thick} if $a,b \geq L$. A \emph{grid of hyperplanes} is the data of two families of hyperplanes $\mathcal{V}= \{V_1, \ldots, V_p \}$ and $\mathcal{H}= \{H_1, \ldots, H_q\}$, such that any $V_i$ is transverse to any $H_j$, and any $V_i$ (resp. $H_j$) separates $V_{i-1}$ and $V_{i+1}$ (resp. $H_{j-1}$ and $H_{j+1}$); such a grid is said \emph{$\delta$-thin} if  $\min ( \# \mathcal{V}, \# \mathcal{H}) \leq \delta$. Notice that an $(n,m)$-grid of hyperplanes is precisely the crossing graph of $[0,n] \times [0,m] \subset \mathbb{R}^2$ endowed with its canonical structure of square complex.

\medskip \noindent
Then, we notice that a cone-off ``killing'' the flat rectangles is hyperbolic:

\begin{thm}\label{killingflatrectangles}
Let $X$ be a CAT(0) cube complex and $Y$ a cone-off of $X$ over a collection of combinatorially convex subcomplexes. If the thick flat rectangles of $X$ are uniformly bounded in $Y$, then $Y$ is hyperbolic.
\end{thm}

Our first application concerns weak relative hyperbolicity. In \cite{MR3217625}, Hagen associates to any CAT(0) cube complex a hyperbolic graph, namely its \emph{contact graph}, and in particular he deduces that cubulable groups are weakly hyperbolic relatively to the hyperplane stabilizers. Up to a quasi-isometry, looking at the contact graph amounts to cone-off the (neighborhood of the) hyperplanes of the cube complex. Using Theorem \ref{killingflatrectangles}, we show that it is possible to choose the hyperplanes we have to cone-off. More precisely, for every $n \geq 0$, we define a class of hyperplanes, the \emph{$n$-combinatorially contracting} hyperplanes, and we prove that the cone-off $\Gamma_nX$ of $X$ over the hyperplanes which are not $n$-combinatorially contracting is hyperbolic. In particular, we are able to slightly improve the result of Hagen:

\begin{cor}
Let $G$ be a group acting geometrically on a CAT(0) cube complex $X$. Then $G$ is weakly hyperbolic relatively to the stabilizers of the non-contracting hyperplanes of $X$.
\end{cor}

\noindent
For instance, applied to the class of right-angled Coxeter groups, we obtain:

\begin{prop}
Let $\Gamma$ be a finite graph. The right-angled Coxeter group $C(\Gamma)$ is weakly hyperbolic relatively to each of the following collection of subgroups
\begin{itemize}
	\item $\{ C(\Gamma_1 \ast \Gamma_2) \mid \Gamma_1 \ast \Gamma_2 \subset \Gamma \ \text{with} \ \Gamma_1,\Gamma_2 \ \text{not complete}\}$,
	\item $\{ \langle \mathrm{star}(u) \rangle \mid u \in \square(\Gamma) \}$,
\end{itemize}
where $\square(\Gamma)$ denotes the set of vertices of $\Gamma$ which belong to an induced square.
\end{prop}

Afterwards, we focus on strong relative hyperbolicity. Therefore, we first prove a criterion to determine when the usual cone-off of a CAT(0) cube complex is fine. Explicitly,

\begin{thm}
Let $X$ be a uniformly locally finite CAT(0) cube complex and $Y$ a usual cone-off of $X$ over a collection $\mathcal{Q}$ of combinatorially convex subcomplexes. If $\mathcal{Q}$ is locally finite (ie., there exist only finitely-many subcomplexes of $\mathcal{Q}$ containing a given edge of $X$) and if there exists a constant $C \geq 0$ such that two subcomplexes of $\mathcal{Q}$ are both intersected by at most $C$ hyperplanes, then $Y$ is fine. Conversely, if $\mathcal{Q}$ is not locally finite or if it contains two subcomplexes both intersected by infinitely many hyperplanes, then $Y$ is not fine. 
\end{thm}

Thus, combining this criterion with Theorem \ref{killingflatrectangles}, we are able to reprove the characterization of the strong relative hyperbolicity of right-angled Coxeter groups stated in \cite[Theorem I]{BHSC}. Given a finite graph $\Gamma$, we find a collection of subgraphs $\mathfrak{J}^{\infty}(\Gamma)$ such that:

\begin{thm}
The right-angled Coxeter group $C(\Gamma)$ is relatively hyperbolic if and only if $\mathfrak{J}^{\infty} (\Gamma) \neq \{ \Gamma \}$. If so, then $C(\Gamma)$ is hyperbolic relatively to $\{ C(\Lambda) \mid \Lambda \in \mathfrak{J}^{\infty} (\Gamma) \}$.
\end{thm}

According to Theorem \ref{hyperbolic}, infinite-dimensional CAT(0) cube complexes cannot be hyperbolic. Thus, a natural cone-off to look at is the cone-off over the high dimensional cubes. In fact, up to a quasi-isometry, it amounts to consider the distance $d_{\infty}$ corresponding to the $\ell_{\infty}$-norm on each cube. Thus, considering the distance $d_{\infty}$ allows us to introduce hyperbolicity in infinite dimensions. In fact, since loosely speaking the cone-off kills the dimension, the criterion suggested by Theorem \ref{hyperbolic} is precisely what we get:

\begin{thm}\label{hyperbolicinfty}
Let $X$ be a CAT(0) cube complex. Then $(X,d_{\infty})$ is hyperbolic if and only if the grid of hyperplanes in $X$ are uniformly thin.
\end{thm}

For instance, infinite-dimensional CAT(0) cube complexes in which the grid of hyperplanes are uniformly thin naturally appear in some infinitely-presented groups which are ``limits'' of hyperbolic groups; for example, they include the cubulations of infinitely-presented $C'(1/6)$ or $C'(1/4)-T(4)$ groups \cite{MR2053602}. In order to exhibit a negatively-curved behaviour of such groups, we determine when the action on this hyperbolic space satisfies an acylindrical condition.

\begin{thm}\label{roughlyacyl}
Let $G$ be a group acting on a complete CAT(0) cube complex $X$. Suppose $(X,d_{\infty})$ hyperbolic. The following statements are equivalent:
\begin{itemize}
	\item[(i)] for every $d \geq 0$, there exists $R \geq 0$ such that, for every vertices $x,y \in X$,
\begin{center}
$d_{\infty}(x,y) \geq R \Rightarrow \# \{ g \in G \mid d_{\infty}(x,gx),d_{\infty}(y,gy) \leq d\} <+ \infty$;
\end{center}
	\item[(ii)] there exists $R \geq 0$ such that, for every vertices $x,y \in X$,
\begin{center}
$d_{\infty}(x,y) \geq R \Rightarrow \# \{ g \in G \mid gx=x, gy=y\} <+ \infty$;
\end{center}
	\item[(iii)] there exists $R \geq 0$ such that, for any hyperplanes $J_1,J_2$ separated by at least $R$ pairwise disjoint hyperplanes, $\mathrm{stab}(J_1) \cap \mathrm{stab}(J_2)$ is finite.
\end{itemize}
\end{thm}

Such an action is useful to prove the acylindrical hyperbolicity of the group (see Section \ref{small} for precise definitions), since any loxodromic isometry of $G$ turns out to be WPD. For instance, we are able to deduce:

\begin{thm}\label{smallcancellation1}
Infinitely-presented $C'(1/4)-T(4)$ groups are acylindrically hyperbolic.
\end{thm}

\begin{thm}\label{smallcancellation2}
Let $G=G_1 \ast \cdots \ast G_n$ be a finitely-generated free product and $R \subset G$ a family satisfying the condition $C'(1/4)-T(4)$. Then the quotient $Q= G / \langle \langle R \rangle \rangle$ is either virtually cyclic or acylindrically hyperbolic.
\end{thm}

The first result may be compared with \cite{C(7)}, where it is proved, in particular, that infinitely-presented C(7) groups are acylindrically hyperbolic.

\medskip
This paper is organised as follows. Section 2 contains the necessary prelimiaries on CAT(0) cube complexes, including combinatorial projections and disc diagrams. Sections 3 and 4 are respectively dedicated to the proofs of Theorems \ref{hyperbolic} and \ref{killingflatrectangles}. In Section 5, we prove the applications toward weak and strong relative hyperbolicities that we mentionned above. Finally, Section 6 contains the proof of \ref{hyperbolicinfty}, Section 7 the proof of Theorem \ref{roughlyacyl}, and these two results are combined in Section 8 to deduce Theorem \ref{smallcancellation1} and Theorem \ref{smallcancellation2}.

\section{Preliminaries}

A \textit{cube complex} is a CW complex constructed by gluing together cubes of arbitrary (finite) dimension by isometries along their faces. Furthermore, it is \textit{nonpositively curved} if the link of any of its vertices is a simplicial \textit{flag} complex (ie., $n+1$ vertices span a $n$-simplex if and only if they are pairwise adjacent), and \textit{CAT(0)} if it is nonpositively curved and simply-connected. See \cite[page 111]{MR1744486} for more information.

Alternatively, CAT(0) cube complexes may be described by their 1-skeletons. Indeed, Chepoi notices in \cite{mediangraphs} that the class of graphs appearing as 1-skeletons of CAT(0) cube complexes coincides with the class of \textit{median graphs}, which we now define.

Let $\Gamma$ be a graph. If $x,y,z \in \Gamma$ are three vertices, a vertex $m$ is called a \textit{median point of $x,y,z$} whenever
\begin{center}
$d(x,y)=d(x,m)+d(m,y)$, $d(x,z)=d(x,m)+d(m,z)$, $d(y,z)=d(y,m)+d(m,z)$.
\end{center}
Notice that, for every geodesics $[x,m]$, $[y,m]$ and $[z,m]$, the concatenations $[x,m] \cup [m,y]$, $[x,m] \cup [m,z]$ and $[y,m] \cup [m,z]$ are also geodesics; furthermore, if $[x,y]$, $[y,z]$ and $[x,z]$ are geodesics, then any vertex of $[x,y] \cap [y,z] \cap [x,z]$ is a median point of $x,y,z$.

The graph $\Gamma$ is \textit{median} if every triple $(x,y,z)$ of pairwise distinct vertices admits a unique median point, denoted by $m(x,y,z)$.

\begin{thm}\emph{\cite[Theorem 6.1]{mediangraphs}} 
A graph is median if and only if it is the 1-skeleton of a CAT(0) cube complex.
\end{thm}

A fundamental feature of cube complexes is the notion of \textit{hyperplane}. Let $X$ be a nonpositively curved cube complex. Formally, a \textit{hyperplane} $J$ is an equivalence class of edges, where two edges $e$ and $f$ are equivalent whenever there exists a sequence of edges $e=e_0,e_1,\ldots, e_{n-1},e_n=f$ where $e_i$ and $e_{i+1}$ are parallel sides of some square in $X$. Notice that a hyperplane is uniquely determined by one of its edges, so if $e \in J$ we say that $J$ is the \textit{hyperplane dual to $e$}. Geometrically, a hyperplane $J$ is rather thought of as the union of the \textit{midcubes} transverse to the edges belonging to $J$.
\begin{center}
\includegraphics[scale=0.4]{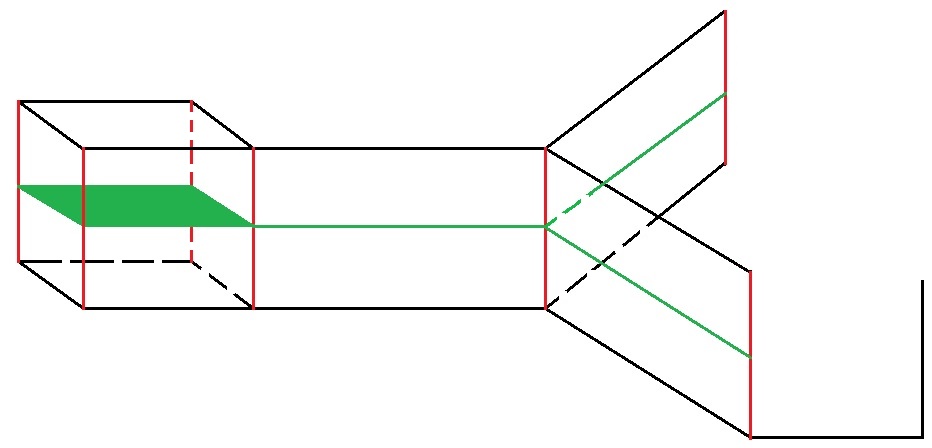}
\end{center}
The \textit{neighborhood} $N(J)$ of a hyperplane $J$ is the smallest subcomplex of $X$ containing $J$, i.e., the union of the cubes intersecting $J$. In the following, $\partial N(J)$ will denote the union of the cubes of $X$ contained in $N(J)$ but not intersecting $J$, and $X \backslash \backslash J= \left( X \backslash N(J) \right) \cup \partial N(J)$. Notice that $N(J)$ and $X \backslash \backslash J$ are subcomplexes of $X$.

\begin{thm}\emph{\cite[Theorem 4.10]{MR1347406}}
Let $X$ be a CAT(0) cube complex and $J$ a hyperplane. Then $X \backslash \backslash J$ has exactly two connected components.
\end{thm}

\noindent
The two connected components of $X \backslash \backslash J$ will be refered to as the \textit{halfspaces} associated to the hyperplane $J$.

\paragraph{Distances $\ell_p$.} 
There exist several natural metrics on a CAT(0) cube complex. For example, for any $p \in (0,+ \infty)$, the $\ell_p$-norm defined on each cube can be extended to a distance defined on the whole complex, the \emph{$\ell_p$-metric}. Usually, the $\ell_1$-metric is refered to as the \emph{combinatorial distance} and the $\ell_2$-metric as the \emph{CAT(0) distance}. Indeed, a CAT(0) cube complex endowed with its CAT(0) distance turns out to be a CAT(0) space \cite[Theorem C.9]{Leary}, and the combinatorial distance between two vertices corresponds to the graph metric associated to the 1-skeleton $X^{(1)}$. In particular, \textit{combinatorial geodesics} are edge-paths of minimal length, and a subcomplex is \textit{combinatorially convex} if it contains any combinatorial geodesic between two of its points.

In fact, the combinatorial metric and the hyperplanes are strongly linked together: the combinatorial distance between two vertices corresponds exactly to the number of hyperplanes separating them \cite[Theorem 2.7]{MR2413337}, and

\begin{thm}\emph{\cite[Corollary 2.16]{MR2413337}}
Let $X$ be a CAT(0) cube complex and $J$ a hyperplane. The two components of $X \backslash \backslash J$ are combinatorially convex, as well as the components of $\partial N(J)$.
\end{thm}

The $\ell_{\infty}$-metric, denoted by $d_{\infty}$, is also of particular interest. Alternatively, given a CAT(0) cube complex $X$, the distance $d_{\infty}$ between two vertices corresponds to the distance associated to the graph obtained from $X^{(1)}$ by adding an edge between two vertices whenever they belong to a common cube. Nevertheless, the distance we obtain stays strongly related to the combinatorial structure of $X$:

\begin{prop}\emph{\cite[Corollary 2.5]{depth}}
Let $X$ be a CAT(0) cube complex and $x,y \in X$ two vertices. Then $d_{\infty}(x,y)$ is the maximal number of pairwise disjoint hyperplanes separating $x$ and $y$.
\end{prop}

\begin{prop}\label{hsgeodesic}
In a CAT(0) cube complex, any half-space is geodesic with respect to the distance $d_{\infty}$.
\end{prop}

\noindent
\textbf{Proof.} Let $X$ be a CAT(0) cube complex, $D$ a half-space and $x,y \in D$ two vertices. To conclude, it is sufficient to prove that a $d_{\infty}$-geodesic between $x$ and $y$ intersecting a minimal number of half-spaces must be included into $D$. Let $\gamma$ be such a geodesic and suppose by contradiction that $\gamma$ is not included into $D$. In particular, there exists a hyperplane $J$ intersected twice by $\gamma$ and a subsegment $\gamma_0 \subset \gamma$ with endpoints in $N(J)$ such that $\gamma_0$ does not intersect twice any hyperplane. By linking two consecutive vertices of $\gamma_0$ by a combinatorial geodesic, we produce a combatorial path $\overline{\gamma_0}$ with the same endpoints as $\gamma_0$; moreover, since $\gamma_0$ does not intersect twice any hyperplane, so does $\overline{\gamma_0}$, ie., $\overline{\gamma_0}$ is a combinatorial geodesic. This implies that 
\begin{center}
$\gamma_0 \subset \overline{\gamma_0} \subset \partial N(J)$, 
\end{center}
by the combinatorial convexity of the components of $\partial N(J)$. Let $\gamma'$ be the $d_{\infty}$-path obtained from $\gamma$ by replacing the subsegment $\gamma_0$ with its image in $N(J)$ by the reflection with respect to $J$. Noticing that $\mathrm{length}(\gamma') \leq \mathrm{length}(\gamma)$, we deduce that $\gamma'$ is a new $d_{\infty}$-geodesic with the same endpoints as $\gamma$. On the other hand, the number of half-spaces intersected by $\gamma'$ is strictly smaller than the same number for $\gamma$, a contradiction. Therefore, $\gamma$ is included into $D$. $\square$

\paragraph{Combinatorial projection}
In CAT(0) spaces, and so in particular in CAT(0) cube complexes with respect to the CAT(0) distance, the existence of a well-defined projection onto a given convex subspace provides a useful tool. Similarly, with respect to the combinatorial distance, it is possible to introduce a \emph{combinatorial projection} onto a combinatorially convex subcomplex, defined by the following result.

\begin{prop}\label{projection}
\emph{\cite[Lemma 1.2.3]{arXiv:1505.02053}} Let $X$ be a CAT(0) cube complex, $C \subset X$ be a combinatorially convex subspace and $x \in X \backslash C$ be a vertex. Then there exists a unique vertex $y \in C$ minimizing the distance to $x$. Moreover, for any vertex of $C$, there exists a combinatorial geodesic from it to $x$ passing through $y$.
\end{prop}

\noindent
For example, we are able to deduce the following result:

\begin{prop}\label{hyperplanseparantcor}
Let $X$ be a CAT(0) cube complex and $C_1,C_2 \subset X$ two combinatorially convex subspaces. If $x \in C_1$ and $y \in C_2$ minimize the distance between $C_1$ and $C_2$ then the hyperplanes separating $x$ and $y$ are precisely those separating $C_1$ and $C_2$.
\end{prop}

\begin{lemma}\label{hyperplan séparant}
Let $X$ be a CAT(0) cube complex and $N \subset X$ a combinatorially convex subspace. Let $p : X \to N$ denote the combinatorial projection onto $N$. Then every hyperplane separating $x$ and $p(x)$ separates $x$ and $N$.
\end{lemma}

\noindent
\textbf{Proof.} According to Proposition \ref{projection}, for all $z \in N$ there exists a combinatorial geodesic $\gamma$ between $x$ and $z$ passing through $p(x)$. Of course, any hyperplane separating $x$ and $p(x)$ meets $\gamma$, so any such hyperplane cannot meet $[p(x),z] \subset \gamma$. Thus, we have proved that no hyperplane separating $x$ and $p(x)$ separates $p(x)$ and some vertex of $N$. $\square$

\medskip \noindent
\textbf{Proof of Proposition \ref{hyperplanseparantcor}.} Clearly, a hyperplane separating $C_1$ and $C_2$ separates $x$ and $y$. Conversely, let $J$ be a hyperplane separating $x$ and $y$, and let $p : X \to C_1$ and $q : X \to C_2$ denote the combinatorial projection onto $C_1$ and $C_2$ respectively. Because $x$ and $y$ minimize the distance between $C_1$ and $C_2$, it follows that $x=p(y)$ and $y=q(x)$. According to Lemma \ref{hyperplan séparant}, we deduce that $J$ separates $x$ and $C_2$, and $y$ and $C_1$. Thus, $J$ separates $C_1$ and $C_2$. $\square$

\medskip \noindent
We conclude this paragraph with this last result:

\begin{prop}\label{proj}
Let $X$ be a CAT(0) cube complex and $C_1,C_2 \subset X$ two combinatorially convex subcomplexes. Let $p : X \to C_2$ denote the combinatorial projection onto $C_2$. Then $p(C_1)$ is a geodesic subcomplex of $C_2$. Moreover, the hyperplanes intersecting $p(C_1)$ are precisely those which intersect both $C_1$ and $C_2$. 
\end{prop}

\begin{lemma}
Let $X$ be a CAT(0) cube complex and $C$ a combinatorially convex subcomplex. Let $p : X \to C$ denote the combinatorial projection onto $C$. For any vertices $x,y \in X$ and any hyperplane $J$ separating $p(x)$ and $p(y)$, $J$ separates $x$ and $y$. In particular, $d(p(x),p(y)) \leq d(x,y)$.
\end{lemma}

\noindent
\textbf{Proof.} Because the hyperplane $J$ separates two vertices of $C$, namely $p(x)$ and $p(y)$, it follows from Lemma \ref{hyperplan séparant} that $J$ cannot separate neither $x$ and $p(x)$ nor $y$ and $p(y)$. Necessarily, $J$ separates $x$ and $y$. $\square$

\medskip \noindent
\textbf{Proof of Proposition \ref{proj}.} Because $d(p(x),p(y)) \leq d(x,y)$, two adjacent vertices in $C_1$ are sent by $p$ onto the same vertex or onto two adjacent vertices. Thus, the image by $p$ of a combinatorial geodesic $\gamma$ in $C_1$ between two vertices $x$ and $y$ defines a combinatorial path $p(\gamma)$ in $C_2$ between $p(x)$ and $p(y)$. Using the previous lemma, we deduce that the length of $p(\gamma)$ is precisely the number of hyperplanes separating $p(x)$ and $p(y)$. Therefore, $p(\gamma)$ is a combinatorial geodesic between $p(x)$ and $p(y)$ included into $p(C_1)$. $\square$

\paragraph{Disc diagrams} 
A fundamental tool to study CAT(0) cube complexes is the theory of \emph{disc diagrams}. For example, they were extensively used by Sageev in \cite{MR1347406} and by Wise in \cite{Wise}. The rest of this section is dedicated to basic definitions and properties of disc diagrams.

\begin{definition}
Let $X$ be a nonpostively curved cube complex. A \emph{disc diagram} is a continuous combinatorial map $D \to X$, where $D$ is a finite contractible square complex with a fixed topological embedding into $\mathbb{S}^2$; notice that $D$ may be \emph{non-degenerated}, ie., homeomorphic to a disc, or may be \emph{degenerated}. In particular, the complement of $D$ in $\mathbb{S}^2$ is a $2$-cell, whose attaching map will be refered to as the \emph{boundary path} $\partial D \to X$ of $D \to X$; it is a combinatorial path. The \emph{area} of $D \to X$, denoted by $\mathrm{Area}(D)$, corresponds to the number of squares of $D$.
\end{definition}

Given a combinatorial closed path $P \to X$, we say that a disc diagram $D \to X$ is \emph{bounded} by $P \to X$ if there exists an isomorphism $P \to \partial D$ such that the following diagram is commutative:
\begin{displaymath}
\xymatrix{ \partial D \ar[rr] & & X \\ P \ar[u] \ar[urr] & & }
\end{displaymath}
According to a classical argument due to van Kampen \cite{vanKampen}, there exists a disc diagram bounded by a given combinatorial closed path if and only if this path is null-homotopic. Thus, if $X$ is a CAT(0) cube complex, then any combinatorial closed path bounds a disc diagram.

\medskip 
As a square complex, a disc diagram contains hyperplanes: they are called \emph{dual curves}. Equivalently, they correspond to the connected components of the reciprocal images of the hyperplanes of $X$. Given a disc diagram $D \to X$, a \emph{nogon} is a dual curve homeomorphic to a circle; a \emph{monogon} is a subpath, of a self-intersecting dual curve, homeomorphic to a circle; an \emph{oscugon} is a subpath of a dual curve whose endpoints are the midpoints of two adjacent edges; a \emph{bigon} is a pair of dual curves intersecting into two different squares.
\begin{figure}[ht!]
\includegraphics[scale=0.4]{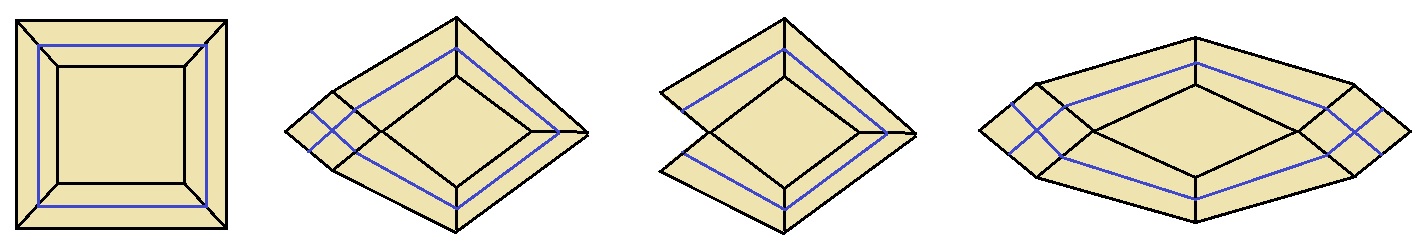}
\caption{From left to right: a nogon, a monogon, an oscugon and a bigon.}
\end{figure}

\begin{thm}\label{Wise}\emph{\cite[Lemma 2.2]{Wise}}
Let $X$ be a nonpositively curved cube complex and $D \to X$ a disc diagram. If $D$ contains a nogon, a monogon, a bigon or an oscugon, then there exists a new disc diagram $D' \to X$ such that:
\begin{itemize}
	\item[(i)] $D'$ is bounded by $\partial D$,
	\item[(ii)] $\mathrm{Area}(D') \leq \mathrm{Aire}(D) -2$.
\end{itemize}
\end{thm}

Let $X$ be a CAT(0) cube complex. A \emph{cycle of subcomplexes} $\mathcal{C}$ is a sequence of subcomplexes $C_1, \ldots, C_r$ such that $C_1 \cap C_r \neq \emptyset$ and $C_i \cap C_{i+1} \neq \emptyset$ for every $1 \leq i \leq r-1$. A disc diagram $D \to X$ is bounded by $\mathcal{C}$ if $\partial D \to X$ can be written as the concatenation of $r$ combinatorial geodesics $P_1, \ldots, P_r \to X$ such that $P_i \subset C_i$ for every $1 \leq i \leq r$. The \emph{complexity} of such a disc diagram is defined by the couple $(\mathrm{Area}(D), \mathrm{length}(\partial D))$, and a disc diagram bounded by $\mathcal{C}$ will be of \emph{minimal complexity} if its complexity is minimal with respect to the lexicographic order among all the possible disc diagrams bounded by $\mathcal{C}$ (allowing modifications of the paths $P_i$).

\medskip 
Our main technical result on disc diagrams is the following, inspired by a result due to Hagen \cite[Lemma 2.11]{MR3217625}. The arguments used here are essentially the same, and mainly come from \cite{Wise}.

\begin{thm}\label{disc diagram}
Let $X$ be a CAT(0) cube complex, $\mathcal{C}=(C_1, \ldots, C_r)$ a cycle of subcomplexes, and $D \to X$ a disc diagram bounded by $\mathcal{C}$. For convenience, write $\partial D$ as the concatenation of $r$ combinatorial geodesics $P_1, \ldots, P_r \to X$ with $P_i \subset C_i$ for every $1 \leq i\leq r$. If the complexity of $D \to X$ is minimal, then:
\begin{itemize}
	\item[(i)] if $C_i$ is combinatorially convex, two dual curves intersecting $P_i$ are disjoint;
	\item[(ii)] if $C_i$ and $C_{i+1}$ are combinatorially convex, no dual curve intersects both $P_i$ and $P_{i+1}$.
\end{itemize}
\end{thm}

\noindent
\textbf{Proof.} Suppose that $C_i$ is combinatorially convex and that there exist two transverse dual curves $c_1$ and $c_2$ intersecting $P_i$. Any dual curve intersecting $P_i$ between $c_1$ and $c_2$ necessarily intersects either $c_1$ or $c_2$, so we may suppose without loss of generality that $c_1$ and $c_2$ intersect $P_i$ along two adjacent edges. If these two edges bound a square inside $D$, then we replace $P_i$ with $P_i'$ as follows: 
\begin{center}
\includegraphics[scale=0.6]{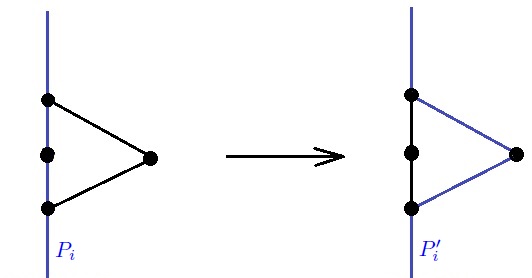}
\end{center}
Because $C_i$ is combinatorially convex, $P_i' \subset C_i$. Thus, we have constructed a new disc diagram $D' \to X$ bounded by $\mathcal{C}$, such that $\mathrm{Area}(D')= \mathrm{Area}(D)-1$, so that $c(D')< c(D)$. This contradicts the minimality of the complexity of $D$. Now, suppose that the edges of $P_i$ dual to $c_1$ and $c_2$ do not bound a square inside $D$. Because the hyperplanes of $X$ do not inter-osculate, the images of these two edges bound a square inside $X$, so that it is possible to construct a new disc diagram $D' \to X$ as follows:
\begin{center}
\includegraphics[scale=0.54]{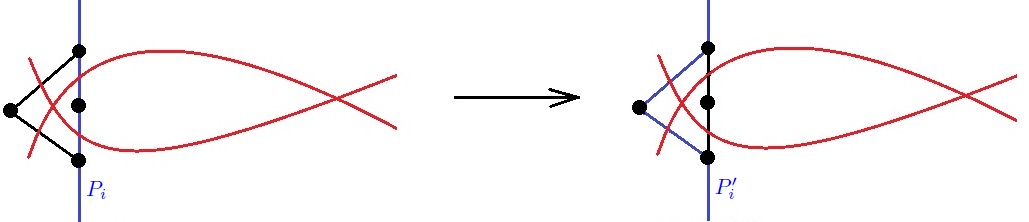}
\end{center}
Similarly, by the combinatorial convexity of $C_i$, we have $P_i' \subset C_i$ so that $D' \to X$ is bounded by $\mathcal{C}$. Noticing that $D'$ contains a bigon, Theorem \ref{Wise} allows us to construct a new disc diagram $D'' \to X$ bounded by $\partial D'$ and satisfying $\mathrm{Area}(D'') \leq \mathrm{Area}(D')-2$. Therefore, we have $\mathrm{Aire}(D'') \leq \mathrm{Aire}(D)-1$, so that $c(D'')<c(D)$ which contradicts the minimality of the complexity of $D$. This concludes the proof of the point $(i)$.

\medskip \noindent
In fact, we have proved a more precise statement which will be useful later:

\begin{fact}\label{fact}
Let $X$ be a CAT(0) cube complex, $\mathcal{C}=(C_1, \ldots, C_r)$ a cycle of subcomplexes, and $D \to X$ a disc diagram bounded by $\mathcal{C}$. For convenience, write $\partial D$ as the concatenation of $r$ combinatorial geodesics $P_1, \ldots, P_r \to X$ with $P_i \subset C_i$ for every $1 \leq i\leq r$. If there exists an index $1 \leq i \leq r$ such that there exist two transverse dual curves intersecting $C_i$, then there exists a new disc diagram $D' \to X$ satisfying $\mathrm{Area}(D') \leq \mathrm{Area}(D)-2$ and whose path boundary is the concatenation of $P_1, \ldots, P_{i-1},P_i',P_{i+1},\ldots, P_r$ where $P_i'$ is a combinatorial geodesic with the same endpoints as $P_i$. 
\end{fact}

\noindent
Now, suppose that $C_i$ and $C_{i+1}$ are combinatorially convex and that there exists a dual curve $c$ intersecting both $P_i$ and $P_{i+1}$; let $e_i$ and $e_{i+1}$ denote respectively the edges of $P_i$ and $P_{i+1}$ dual to $c$. Because two dual curves intersecting either $P_i$ or $P_{i+1}$ are necessarily disjoint according to the point $(i)$, no dual curve intersecting $P_i \cup P_{i+1}$ between $e_i$ and $e_{i+1}$ can intersect $c$, so we may suppose without loss of generality that $e_i$ and $e_{i+1}$ are adjacent. If $e_i \neq e_{i+1}$, then $D$ contains an oscugon so that it is possible to decrease the area of $D$ without disturbing its boundary path according to Theorem \ref{Wise}, contradicting the minimality of the complexity of $D$. Thus, $e_i = e_{i+1}$. Now, replace $P_i$ and $P_{i+1}$ respectively with $P_i'$ and $P_{i+1}'$ as follows:
\begin{center}
\includegraphics[scale=0.6]{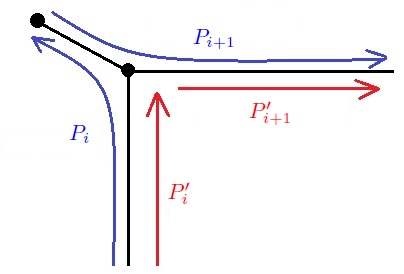}
\end{center}
We have constructed a new disc diagram $D' \to X$ bounded by $\mathcal{C}$ and satisfying $\mathrm{Area}(D')= \mathrm{Area}(D)$, $\mathrm{length}(\partial D') = \mathrm{length}(\partial D)-1$. In particular, we have $c(D')<c(D)$, contradicting the minimality of the complexity of $D$. $\square$

\medskip \noindent
We conclude the section by a last result on disc diagrams which will be useful later.

\begin{prop}\label{disc embedding}
Let $X$ be a CAT(0) cube complex and $D \to X$ a disc diagram which does not contain any bigon. With respect to the combinatorial metrics, $\varphi : D \to X$ is an isometric embedding if and only if every hyperplane of $X$ induces at most one dual curve of $D$. 
\end{prop}

\begin{lemma}
Let $D \to X$ be a disc diagram which does not contain any bigon and $x,y \in D$ two vertices. The combinatorial distance between $x$ and $y$ in $D$ is equal to the number of dual curves separating $x$ and $y$. 
\end{lemma}

\noindent
\textbf{Proof.} Let $c(x,y)$ denote the number of dual curves separating $x$ and $y$. Of course, any combinatorial path between $x$ and $y$ must intersect each of the $c(x,y)$ dual curves separating $x$ and $y$, hence $d(x,y) \geq c(x,y)$. Conversely, to prove that $d(x,y) \leq c(x,y)$, it is sufficient to show that any combinatorial path of minimal length between $x$ and $y$ intersects each dual curve at most once. Let us consider a combinatorial path $\gamma$ between $x$ and $y$ which intersects a dual curve at least twice; in particular, there exist two edges $e$ and $f$ of this path dual to the same dual curve. We choose this dual curve innermost, so that, if $\gamma_2$ denotes the subpath of $\gamma$ between $e$ and $f$, then no dual curve intersects $\gamma_2$ twice. Let $\gamma_0,\gamma_1$ denote the two combinatorial paths fellow-traveling our dual curve between the endpoints of $e$ and $f$.
\begin{center}
\includegraphics[scale=0.6]{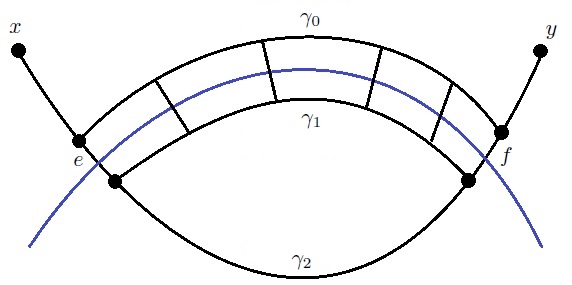}
\end{center}
In particular, any dual curve intersecting $\gamma_2$ intersects $\gamma_1$. Conversely, because $D$ does not contain any bigon, any dual curve intersecting $\gamma_1$ intersects necessarily $\gamma_2$. Consequently, $\mathrm{length}(\gamma_1)= \mathrm{length}(\gamma_2)$. By replacing the subpath $e \cup \gamma_2 \cup f$ with $\gamma_0$ in our path, we get a new combinatorial path whose length is smaller, since
\begin{center}
$\mathrm{length}(e \cup \gamma_2 \cup f)= \mathrm{length}(\gamma_2)+2 = \mathrm{length}(\gamma_1)+2= \mathrm{length}(\gamma_0)+2 > \mathrm{length}(\gamma_0)$.
\end{center}
We have proved that any combinatorial path between $x$ and $y$ which intersects a dual curve twice may be shortened. This concludes the proof. $\square$

\medskip \noindent
\textbf{Proof of Proposition \ref{disc embedding}.} Suppose that there exist two dual curves $c_1$ and $c_2$ induced by the same hyperplane of $X$. Let $x,y \in D$ be two vertices separated by $c_1$ and $c_2$. The combinatorial distance between $\varphi(x)$ and $\varphi(y)$ corresponds to the number of hyperplanes separating them; each of these hyperplanes induces a dual curve in $D$ separating $x$ and $y$, so there exist at least $d(\varphi(x),\varphi(y))$ dual curves separating $x$ and $y$, hence $d(x,y) \geq d(\varphi(x),\varphi(y))$. On the other hand, we know that at least two of these dual curves, namely $c_1$ and $c_2$, are induced by the same hyperplane, hence $d(x,y) > d(\varphi(x),\varphi(y))$. Consequently, $\varphi : D \to X$ is not an isometric embedding.

\medskip \noindent
Conversely, suppose that each hyperplane of $X$ induces at most one dual curve on $D$. Now fix two vertices $x,y \in D$ together with a combinatorial path $\gamma$ of minimal length between $x$ and $y$ in $D$. According to the previous lemma, any dual curve intersects $\gamma$ at most once. Thus, the only possibility for a hyperplane of $X$ to intersect $\varphi(\gamma)$ twice is to induce two different dual curve intersecting $\gamma$, which is impossible by our hypothesis. We conclude that $\varphi (\gamma)$ is a combinatorial geodesic in $X$, ie., $d(\varphi(x), \varphi(y)) \geq d(x,y)$. We have already noticed in the previous paragraph the inequality $d(x,y) \geq d(\varphi(x),\varphi(y))$, so $d(x,y)=d( \varphi(x),\varphi(y))$. Consequently, $\varphi : D \to X$ is an isometric embedding. $\square$

\begin{cor}\label{flat rectangle}
Let $X$ be a CAT(0) cube complex and $\mathcal{C}$ a cycle of four combinatorially convex subcomplexes. If $D \to X$ is a disc diagram of minimal complexity bounded by $\mathcal{C}$, then $D$ is combinatorially isometric to a rectangle $[0,a] \times [0,b] \subset \mathbb{R}^2$ and $D \to X$ is an isometric embedding.
\end{cor}

\noindent
\textbf{Proof.} By definition, the boundary path $\partial D \to X$ can be written as the concatenation of four combinatorial geodesics $P_1,P_2,P_3,P_4 \to X$. It follows directly from Theorem \ref{disc diagram} that two dual curves intersecting $P_i$ are necessarily disjoint and that any dual curve intersecting $P_i$ necessarily intersects $P_{i+2}$ (taking $i$ modulo $4$). In particular, $\mathrm{length}(P_1)= \mathrm{length}(P_3)$ and $\mathrm{length}(P_2)=\mathrm{length}(P_4)$; let us say that a dual curve intersecting $P_1$ (resp. $P_2$) is \emph{vertical} (resp. \emph{horizontal}) and let $x_0$ denote the common endpoint of $P_1$ and $P_4$. Now, the map $x \mapsto (v(x),h(x))$, where $v(x)$ (resp. $h(x)$) denotes the number of vertical (resp. horizontal) dual curves separating $x_0$ and $x$, defines a combinatorial isometric embedding $D \hookrightarrow \mathbb{Z}^2$ and its image is a rectangle. Finally, it follows easily from Proposition \ref{disc embedding} that $D \to X$ is an isometric embedding: two dual curves in $D$ either intersects the same $P_i$, so they are necessarily induced by two different hyperplanes since $P_i$ is a combinatorial geodesic, or they are transverse, so that they map onto two transverse hyperplanes, which are necessarily distinct since a CAT(0) cube complex does not contain self-intersecting hyperplanes. $\square$

\section{Hyperbolicity with respect to the \texorpdfstring{$\ell_1$}{l1} and \texorpdfstring{$\ell_2$}{l2} metrics}

\noindent
In this section, we are interested in determining when a CAT(0) cube complex is hyperbolic with respect to the combinatorial or CAT(0) distance. It is worth noticing that an $n$-cube contains a geodesic triangle which is not $(n-1)$-thin with respect to the combinatorial distance, and a geodesic triangle which is not $(\sqrt{n}/2-1)$-thin with respect to the CAT(0) metric; therefore, an infinite-dimensional CAT(0) cube complex is hyperbolic neither with respect to the combinatorial distance nor with respect to the CAT(0) distance. Thus, the hyperbolicity happens only in finite dimension, where these two distances are quasi-isometric, so that being hyperbolic does not depend on the metric we choose. For convenience, we will use the combinatorial distance.

\begin{definition}
A \emph{grid of hyperplanes} is the data of two families of hyperplanes $\mathcal{V}= \{V_1, \ldots, V_p \}$ and $\mathcal{H}= \{H_1, \ldots, H_q\}$, such that any $V_i$ is transverse to any $H_j$, and any $V_i$ (resp. $H_j$) separates $V_{i-1}$ and $V_{i+1}$ (resp. $H_{j-1}$ and $H_{j+1}$); such a grid is said \emph{$\delta$-thin} if  $\min ( \# \mathcal{V}, \# \mathcal{H}) \leq \delta$.
\end{definition}

\begin{definition}
Let $X$ be a CAT(0) cube complex. A \emph{flat rectangle} is a combinatorially geodesic subcomplex of $X$ isomorphic to some rectangle $[0,a] \times [0,b] \subset \mathbb{R}^2$, with $a,b \geq 1$, endowed with its canonical structure of square complex; it is \emph{$L$-thick} if $a,b \geq L$ and \emph{$L$-thin} if $\min(a,b) \leq L$.
\end{definition}

\noindent
Our main criterion is the following:

\begin{thm}\label{critère d'hyperbolicité}
Let $X$ be a CAT(0) cube complex. The following are equivalent:
\begin{itemize}
	\item[(i)] $X$ is hyperbolic;
	\item[(ii)] the flat rectangles in $X$ are uniformly thin;
	\item[(iii)] $X$ is finite-dimensional and the grids of hyperplanes of $X$ are uniformly thin.
\end{itemize}
\end{thm}

\begin{remark}
The equivalence $(i) \Leftrightarrow (iii)$ may be compared with \cite[Theorem 7.3]{MR3217625}, where Hagen proves that a cocompact CAT(0) cube complex is hyperbolic if and only if its crossing graph does not contain a complete bipartite subgraph $K_{\infty,\infty}$. Recall that the crossing graph of a CAT(0) cube complex $X$ is the graph whose vertices are the hyperplanes of $X$ and whose edges link two transverse hyperplanes. 
\end{remark}

\noindent
In fact, a simple criterion of hyperbolicity for median graphs has already been proved by Sigarreta, expressed in terms of bigons \cite{MR3146600}. (Notice that Papasoglu proved that the same statement holds for every graph \cite{Papasoglu}, but the argument of Sigarreta, reproduced below for completeness, is completely elementary.) A fortiori, this produces a criterion of hyperbolicity for CAT(0) cube complexes, which will be used to prove the implication $(iii) \Rightarrow (i)$ of Theorem \ref{critère d'hyperbolicité}.

\begin{definition}
A \textit{bigon} is a pair of combinatorial geodesics $\gamma_1,\gamma_2$ with the same endpoints. It is \textit{$\delta$-thin} if the Hausdorff distance between $\gamma_1$ and $\gamma_2$ is at most $\delta$.
\end{definition}

\begin{lemma}\label{bigone}
\emph{\cite[Theorem 2.5]{MR3146600}} A CAT(0) cube complex is hyperbolic if and only if its bigons are uniformly thin.
\end{lemma}

\noindent
\textbf{Proof.} It is known that, in a $\delta$-hyperbolic geodesic space, the Hausdorff distance between two geodesics with the same endpoints is bounded by a constant $C(\delta)$ depending only on $\delta$. Therefore, the bigons in a hyperbolic graph are uniformly thin. 

\medskip \noindent
Conversely, suppose there exists some $\delta>0$ such that the bigons of $X$ are all $\delta$-thin. Let $[x,y,z]$ be a geodesic triangle. Let $m$ denote the median point of $\{x,y,z\}$. Let $p \in [x,y,z]$ be a point on our triangle, say $p \in [x,y]$. Because the bigon $\{[x,y] , [x,m] \cup [m,y]\}$ is $\delta$-thin, there exists a point $p' \in [x,m] \cup [m,y]$, say $p' \in [x,m]$, such that $d(p,p') \leq \delta$. Again, because the bigon $\{[x,m] \cup [m,z]\cup [x,z]\}$ is $\delta$-thin, there exists a point $p'' \in [x,z]$ such that $d(p',p'') \leq \delta$. Thus,
\begin{center}
$d(p,[x,z]) \leq d(p,p'') \leq d(p,p')+d(p',p'') \leq 2\delta$.
\end{center}
Therefore, the triangle $[x,y,z]$ is $2\delta$-thin. We conclude that $X$ is $2\delta$-hyperbolic. $\square$

\medskip \noindent
Let $\mathrm{Ram}(d)$ denote the Ramsey number defined by: if $K$ is a complete graph with at least $\mathrm{Ram}(d)$ vertices, then any label of the vertices of $K$ with two colors contains a unicolor subgraph containing at least $d+1$ vertices. A simple and well-known application of Ramsey numbers in the context of CAT(0) cube complexes is: 

\begin{lemma}\label{Ramsey}
Let $X$ be a CAT(0) cube complex and let $\mathcal{H}$ be a collection of hyperplanes which does not contain $d+1$ pairwise transverse hyperplanes. (For instance, this happens when $\dim X \leq d$.) If $\mathcal{H}$ has at least $\mathrm{Ram}(d)$ hyperplanes, then its contains a family of at least $d+1$ pairwise disjoint hyperplanes.
\end{lemma}

\noindent
\textbf{Proof.} Let $\mathcal{H}$ be a collection of at least $\mathrm{Ram}(d)$ hyperplanes. Let $K$ be the complete graph whose set of vertices is $\mathcal{H}$ and whose edges are red (resp. blue) if they link two transverse (resp. disjoint) hyperplanes. By definition of $\mathrm{Ram}(d)$, $K$ contains a unicolor subgraph with at least $d+1$ vertices, ie., $\mathcal{H}$ has to contain at least $d+1$ hyperplanes which are either pairwise transverse or pairwise disjoint. But the first possibility is excluded by our hypotheses, so that our lemma follows. $\square$

\medskip \noindent
\textbf{Proof of Theorem \ref{critère d'hyperbolicité}.} The implication $(i) \Rightarrow (ii)$ is clear.

\medskip \noindent
Let us prove $(ii) \Rightarrow (iii)$. Suppose that all flat rectangles in $X$ are $L$-thin and let $(\mathcal{V}=\{V_1, \ldots, V_n\},\mathcal{H}=\{H_1, \ldots, H_m\})$ be an $(n,m)$-grid of hyperplanes. In particular, $\mathcal{C}=(N(V_1),N(H_1),N(V_n),N(H_m))$ defines a cycle of subcomplexes. Let $D \to X$ be a disc diagram of minimal complexity bounded by $\mathcal{C}$. According to Corollary \ref{flat rectangle}, we find a combinatorial isometric embedding $[0,a] \times [0,b] \hookrightarrow X$. Furthermore, because $V_2, \ldots, V_{n-1}$ (resp. $H_{1},\ldots,H_{m-1}$) separate $V_1$ and $V_n$ (resp. $H_1$ and $H_m$), we deduce that $a \geq n-1$ (resp. $b \geq m-1$). We conclude that 
\begin{center}
$\min(n,m) \leq 1+ \min(a,b) \leq L+1$. 
\end{center}
Thus, the grid of hyperplanes in $X$ are uniformly thin. Finally, to prove that $X$ is necessarily finite-dimensional, it is sufficient to notice that the cube $[0,1]^n$ contains a flat rectangle which is not $\lfloor n/2 \rfloor$-thin. For instance, if
\begin{center}
$e_p(i)= ( \underset{i \ \text{coordinates} }{ \underbrace{1,\ldots, 1} }, \underset{ p-i \ \text{coordinates} }{ \underbrace{0, \ldots, 0} })$,
\end{center}
then the vertices $( e_{ \lfloor n/2 \rfloor } (i), e_{ n- \lfloor n/2 \rfloor } (j) ) \in [0,1]^n$, where $0 \leq i \leq \lfloor n/2 \rfloor$ and $0 \leq j \leq n- \lfloor n/2 \rfloor$, span a flat rectangle $[0, \lfloor n/2 \rfloor ] \times [ 0, n- \lfloor n/2 \rfloor ] \hookrightarrow [0,1]^n$.

\medskip \noindent
Finally, we prove $(iii) \Rightarrow (i)$. Suppose that $\dim X <+ \infty$ and that the grid of hyperplanes in $X$ are all $C$-thin for some $C \geq 1$; without loss of generality, we may suppose that $C \geq \dim X$. According to Lemma \ref{bigone}, it is sufficient to prove that the bigons of $X$ are uniformly thin in order to deduce that $X$ is hyperbolic. So let $\gamma_1,\gamma_2$ be two combinatorial geodesics with the same endpoints $x,y$. If $a \in \gamma_1$, let $b \in \gamma_2$ denote the vertex of $\gamma_2$ satisfying $d(x,a)=d(x,b)$. Let $\mathcal{H}_1$ denote the set of hyperplanes meeting $[x,a] \subset \gamma_1$ and separating $a$ and $b$, $\mathcal{H}_2$ the set of hyperplanes meeting $[x,b] \subset \gamma_2$ and separating $a$ and $b$, and finally $\mathcal{H}_3$ the set of hyperplanes separating $x$ and $\{a,b\}$. 
\begin{center}
\includegraphics[scale=0.6]{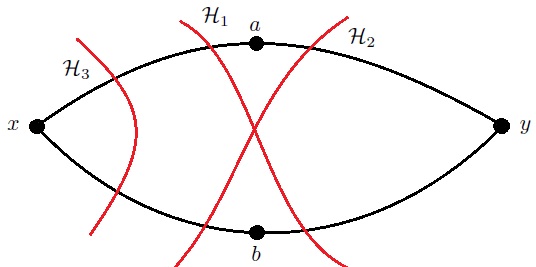}
\end{center}
Because the hyperplanes of $\mathcal{H}_1$ separate $b$ and $\{a,y\}$, and the hyperplanes of $\mathcal{H}_2$ separate $a$ and $\{b,y\}$, we deduce that any hyperplane of $\mathcal{H}_1$ is transverse to any hyperplane of $\mathcal{H}_2$. Now, notice that
\begin{center}
$\# \mathcal{H}_3+ \# \mathcal{H}_1= d(x,a)=d(x,b)= \# \mathcal{H}_3+ \# \mathcal{H}_2$,
\end{center}
whence $\# \mathcal{H}_1 = \# \mathcal{H}_2$. Let $p$ denote this common cardinality. If $p \geq \mathrm{Ram}(C)$, then Lemma \ref{Ramsey} implies that $\mathcal{H}_1$ and $\mathcal{H}_2$ each contain a subfamily of at least $C+1$ pairwise disjoint hyperplanes, producing a $(C+1,C+1)$-grid of hyperplanes and thus contradicting our hypotheses. Therefore, $p \leq \mathrm{Ram}(C)$. We conclude that
\begin{center}
$d(a,b)= \# \mathcal{H}_1 + \# \mathcal{H}_2 = 2p \leq 2\mathrm{Ram}(C)$.
\end{center}
We have proved that the bigons of $X$ are all $2\mathrm{Ram}(C)$-thin, so that $X$ is hyperbolic (in fact, $4\mathrm{Ram}(C)$-hyperbolic). $\square$

\medskip \noindent
In particular, an immediate corollary of Theorem \ref{critère d'hyperbolicité} is that a CAT(0) cube complex whose transversality graph is uniformly locally finite is hyperbolic. In fact, we can even prove more:

\begin{prop}\label{quasi-tree}
A CAT(0) cube complex whose transversality graph is uniformly locally finite is quasi-isometric to a tree.
\end{prop}

\noindent
\textbf{Proof.} Let $X$ be a CAT(0) cube complex whose transversality graph is uniformly locally finite, ie., there exists a constant $N$ such that any hyperplane of $X$ is transverse to at most $N$ other hyperplanes. This implies that the diameter of any hyperplane $J$ is at most $N$. Indeed, if $x,y$ are two vertices which belongs to the same connected component of $\partial N(J)$, then any hyperplane separating $x$ and $y$ is transverse to $J$, hence $d(x,y) \leq N$. Now, it is easy to apply the following criterion:

\medskip \noindent
\textbf{Bottleneck criterion.} \cite{bottleneck} \emph{A geodesic metric space $(S,d)$ is quasi-isometric to a tree if and only if there exists a constant $\delta>0$ such that, any pair of points $x,y \in S$ admits a midpoint $m=m(x,y)$ with the property that any continuous path between $x$ and $y$ intersects the ball $B(m,\delta)$.}

\medskip \noindent
Indeed, if $x,y \in X$ are two points then fix any midpoint $m=m(x,y)$ and let $J$ denote a hyperplane separating $x$ and $y$ such that $m \in N(J)$. Now, let $\gamma$ be any continous path between $x$ and $y$. Of course, $\gamma$ has to pass through $J$, so that it finally intersects the ball $B(m,N) \supset J$. $\square$

\begin{remark}
When $X$ is a CAT(0) cube complex obtained by cubulating a codimension-one subgroup, Niblo characterized the (uniform) local finiteness of the transversality graph of $X$ from an algebraic viewpoint thanks to its \emph{splitting obstruction} \cite[Proposition 7]{SplittingObstruction}, and used it to establish a splitting result in the spirit of Stallings Theorem \cite[Theorem B]{SplittingObstruction}.
\end{remark}

\begin{remark}
The hypothesis of uniform local finiteness in the statement of Proposition \ref{quasi-tree} cannot be weakened into a hypothesis of local finitess. Indeed, let $X$ be a subcomplex of $\mathbb{R}^2$ delimited by two combinatorial rays $r=(r_x,r_y)$ and $\rho=(\rho_x,\rho_y)$ with
\begin{itemize}
	\item $r_x,r_y,\rho_x,\rho_y : \mathbb{N} \to \mathbb{N}$,
	\item $r_x(0)=r_y(0)= \rho_x(0)=\rho_y(0)=0$,
	\item $r_x(t), r_y(t),\rho_x(t),\rho_y(t) \underset{t \to + \infty}{\longrightarrow} + \infty$,
	\item the Hausdorff distance between $r$ and $\rho$ is infinite.
\end{itemize}
Then the transversality graph of $X$ is locally finite but $X$ is not hyperbolic. 
\end{remark}

\section{Killing flat rectangles}

\noindent
In this section, we prove that a cone-off of a CAT(0) cube complex ``killing'' the flat rectangles of Theorem \ref{critère d'hyperbolicité} is hyperbolic. This is the main result of this section.

\begin{thm}\label{KFR}
Let $X$ be a CAT(0) cube complex and $Y$ a cone-off of $X$ over a collection of combinatorially convex subcomplexes. If there exist some constants $C,L \geq 1$ such that the diameter of any $L$-thick flat rectangle of $X$ is at most $C$ in $Y$, then $Y$ is $\delta$-hyperbolic for some $\delta$ depending only on $C$ and $L$.
\end{thm}

\noindent
To prove the hyperbolicity, we will use the following criterion:

\begin{lemma}\label{thinbigon}
Let $X$ be a CAT(0) cube complex and $Y$ a cone-off of $X$ over a collection of combinatorially convex subcomplexes. If there exists some $C \geq 1$ such that the bigons of $X$ are all $C$-thin in $Y$, then $Y$ is $\delta$-hyperbolic for some $\delta$ depending only on $C$.
\end{lemma}

\noindent
\textbf{Proof.} For any vertices $x,y \in Y$, let $\eta(x,y) \subset Y$ denote the union of all the combinatorial geodesics between $x$ and $y$ in $X$. Now, we want to apply the following criterion, due to Bowditch \cite[Proposition 3.1]{Bowditchcriterion}:

\begin{prop}
Let $T$ be a graph and $D \geq 0$. Suppose that a connected subgraph $\eta(x, y)$, containing $x$ and $y$, is associated to any pair of vertices $(x, y) \in T^2$ such that:
\begin{itemize}
	\item for any vertices $x,y \in T$, $d(x,y) \leq 1$ implies $\mathrm{diam}~ \eta(x,y) \leq D$;
	\item for any vertices $x,y,z \in T$, we have $\eta(x,y) \subset ( \eta(x,z) \cup \eta(z,y) )^{+D}$.
\end{itemize}
Then $T$ is $\delta$-hyperbolic for some $\delta$ depending only on $D$.
\end{prop}

\noindent
The first condition is immediate: if $d_Y(x,y) \leq 1$, either $x=y$ and $\eta(x,y)$ is a vertex; or $x$ and $y$ are linked by an edge in $X$, so that they are the only vertices of $\eta(x,y)$; or $x$ and $y$ are linked by an edge in $Y$ which does not belong to $X$, so that there exists a combinatorially convex subcomplex $Q \subset X$ containing $x$ and $y$, and a fortiori $\eta(x,y)$, whose vertices are pairwise linked by an edge in $Y$. Therefore, $d_Y(x,y) \leq 1$ implies $\mathrm{diam}_Y ~ \eta(x,y) \leq 1$. 

\medskip \noindent
Let $x,y,z \in Y$ be three vertices. The second condition of Bowditch's criterion is trivially satisfied if $x,y,z$ are not pairwise distinct. Thus, we may suppose without loss of generality that they are pairwise distinct. Let $m=m(x,y,z)$ denote their median vertex, and fix three combinatorial geodesics $[m,x]$, $[m,y]$ and $[m,z]$ in $X$. Let $w \in \eta(x,y)$, ie., $w$ belongs to some combinatorial geodesic $[x,y]$ between $x$ and $y$ in $X$. By our hypothesis, the bigon $\{ [x,m] \cup [m,y], [x,y] \}$ is $C$-thin in $Y$, so there exists a vertex $w' \in [x,m] \cup [m,y]$ such that $d_Y(w,w') \leq C$. Without loss of generality, say that $w' \in [x,m]$; otherwise, just switch the names of $x$ and $y$. Similarly, the bigon $\{ [x,m] \cup [m,z], [x,z]\}$ is $C$-thin in $Y$, so there exists a vertex $w'' \in [x,z] \subset \eta(x,z)$ such that $d_Y(w',w'') \leq C$. Thus,
\begin{center}
$d_Y(w, \eta(x,z)) \leq d_Y(w,w'') \leq d_Y(w,w')+d_Y(w',w'') \leq 2C$.
\end{center}
We have proved that $\eta(x,y) \subset ( \eta(x,z) \cup \eta(z,y))^{+2C}$. This concludes the proof. $\square$

\medskip \noindent
\textbf{Proof of Theorem \ref{KFR}.} According to the previous lemma, it is sufficient to prove that any bigon of $X$ is $\max(2L,C)$-thin in $Y$. 

\medskip \noindent
Let $x,y \in Y$ be two vertices and $\gamma_1,\gamma_2$ two combinatorial geodesics between $x$ and $y$ in $X$. Fix some vertex $z \in \gamma_1$ and let $z' \in \gamma_2$ be the vertex of $\gamma_2$ defined by $d_X(x,z')=d_X(x,z)$. We claim that $d_Y(z,z') \leq \max(2L,C)$.
\begin{center}
\includegraphics[scale=0.6]{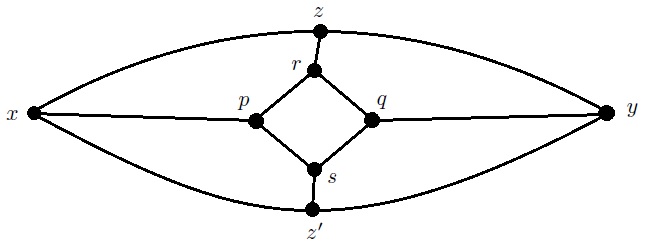}
\end{center}
Let $p$ (resp. $q$) denote the median vertex of $\{x,z,z'\}$ (resp. $\{y,z,z'\}$), and let $r$ (resp. $s$) be the median vertex of $\{ z,p,q\}$ (resp. $\{ z',p,q \}$). Notice that
\begin{center}
$d_X(x,z)=d_X(x,p)+d_X(p,z)=d_X(x,p)+d_X(p,r)+d_X(r,z)$. 
\end{center}
Therefore, if we fix three combinatorial geodesics $[x,p]$, $[p,r]$, $[r,z]$, and if $[x,z]$ denotes the subsegment of $\gamma_1$ between $x$ and $z$, then $\{[x,p] \cup [p,r] \cup [r,z], [x,z] \}$ is a bigon in $X$. Thus, if $J$ is a hyperplane separating $z$ and $r$, necessarily $J$ intersects $[x,z] \subset \gamma_1$. The same argument with respect to $y,r,q,z$ shows that $J$ must intersect $[z,y] \subset \gamma_1$. In particular, $J$ intersects $\gamma_1$ twice: this is impossible since $\gamma_1$ is a combinatorial geodesic. We conclude that no hyperplane separates $z$ and $r$, ie., $r=z$. Similarly, we prove that $s=z'$. In fact, we have proved:

\begin{fact}\label{median}
For any bigon $\{\gamma_1,\gamma_2\}$ with endpoints $x,y$ and for any $z \in \gamma_1$, $z' \in \gamma_2$, we have
\begin{center}
$m(z,m(x,z,z'),m(y,z,z'))=z$.
\end{center}
\end{fact}

\noindent
Consequently, for any choices of four combinatorial geodesics $[p,z]$, $[z,q]$, $[q,z']$ and $[z',p]$, the concatenations $[p,z] \cup [z,q]$, $[z,q] \cup [q,z']$, $[q,z'] \cup [z',p]$ and $[z',p] \cup [p,z]$ are combinatorial geodesics. 

\medskip \noindent
Let $D \to X$ be a disc diagram bounded by a $4$-gon $P=P(p,z,q,z')$, ie., bounded by four combinatorial geodesics $[p,z]$, $[z,q]$, $[q,z']$ and $[z',p]$; suppose that its complexity is minimal among all the disc diagrams bounded by a similar $4$-gon. According to Fact \ref{fact}, if there exist two transverse dual curves in $D$ both intersecting the same side of $P$, then it is possible to find a new disc diagram $D' \to X$ bounded by another $4$-gon $P'=P'(p,z,q,z')$ such that $\mathrm{Area}(D')< \mathrm{Area}(D)$. Thus, two dual curves in $D$ intersecting the same side of $P$ are necessarily disjoint. Furthermore, because the concatenations $[p,z] \cup [z,q]$, $[z,q] \cup [q,z']$, $[q,z'] \cup [z',p]$ and $[z',p] \cup [p,z]$ are combinatorial geodesics, no dual curve in $D$ intersects two adjacent sides of $P$. Consequently, $D$ is a rectangle. In particular, two dual curves in $D$ either are transverse or intersect the same side of $P$, so that a hyperplane of $X$ cannot induce two dual curves in $D$: we deduce from Proposition \ref{disc embedding} that $D \to X$ is an isometric embedding.

\medskip \noindent
Let $\ell=d_X(p,z)=d_X(q,z')$. Notice that
\begin{center}
$d_X(p,z')=d_X(x,z')-d_X(x,p)=d_X(x,z)-d_X(x,p)=d_X(p,z)=\ell$,
\end{center}
and similarly $d_X(q,z)=\ell$. Thus, $D$ defines an $\ell$-thick flat rectangle in $X$. Two cases may happen: either $\ell< L$ and $d_Y(z,z') \leq d_X(z,z')=2 \ell < 2L$; or $\ell \geq L$ and $d_Y(z,z') \leq C$. Therefore, $d_Y(z,z') \leq \max(2L,C)$. $\square$

\section{Relative hyperbolicity}

\subsection{Weak relative hyperbolicity}

\noindent
In the previous section, we have seen how a CAT(0) cube complex becomes hyperbolic when some of its subcomplexes are coned-off. In this section, we apply this criterion to cone-off hyperplanes and finally deduce a weak relative hyperbolicity of cubulable groups. 

\begin{definition}
Let $X$ be a CAT(0) cube complex. A hyperplane $J$ is \emph{$n$-combinatorially contracting} if $\dim J < n$ and $J$ does not belong to an $(n,n)$-grid of hyperplanes in $X$.
\end{definition}

\begin{definition}
Let $X$ be a CAT(0) cube complex and $n \geq 0$ an integer. The \emph{$n$-th contracting graph} of $X$, denobed by $\Gamma_n X$, is the cone-off of $X$ over the neighborhoods of the hyperplanes which are not $n$-combinatorially contracting. \\
If $X$ has finitely-many orbits of hyperplanes under the action of $\mathrm{Aut}(X)$, then the sequence $(\Gamma_nX)$ is eventually constant to some graph $\Gamma_{\infty} X$, which will be called the \emph{contracting graph} of $X$. 
\end{definition}

\begin{prop}\label{contracting}
Let $X$ be a CAT(0) cube complex. Its $n$-th contracting graph $\Gamma_n X$ is $\delta(n)$-hyperbolic, where $\delta(n)$ depends only on $n$.
\end{prop}

\noindent
\textbf{Proof.} According to Theorem \ref{KFR}, it is sufficient to prove that any $\mathrm{Ram}(n)$-thick flat rectangle in $X$ has diameter at most $4 \mathrm{Ram}(n)+3$ in $\Gamma_n X$. For convenience, let $d_n$ denote the distance in $\Gamma_n X$.

\medskip \noindent
Let $Q$ be a $\mathrm{Ram}(n)$-thick flat rectangle in $X$. The set of the hyperplanes intersecting $Q$ can be naturally written as a disjoint union $\mathcal{H} \sqcup \mathcal{V}$ where two hyperplanes both in $\mathcal{H}$ or $\mathcal{V}$ are disjoint in $Q$. 

\medskip \noindent
\underline{Case 1:} $\mathcal{H}$ or $\mathcal{V}$ contains $n$ pairwise transverse hyperplanes. Say that $\mathcal{H}$ contains $n$ pairwise transverse hyperplanes $H_1, \ldots,H_n$. A fortiori, for any hyperplane $V \in \mathcal{V}$, $\{V,H_1, \ldots, H_n\}$ is a collection of $n+1$ pairwise transverse hyperplanes, hence $\dim V \geq n$; therefore, no hyperplane of $\mathcal{V}$ is $n$-combinatorially contracting. If $x_1,x_2 \in Q$ are two vertices, let $V_1,V_2 \in \mathcal{V}$ be two hyperplanes adjacent to $x_1,x_2$ respectively, and let $y_1 \in N(V_1) \cap N(H_1)$, $y_2 \in N(V_2) \cap N(H_1)$ be two vertices. Because $V_1,V_2,H_1$ are not $n$-combinatorially contracting, we deduce that
\begin{center}
$d_n(x_1,x_2) \leq d_n(x_1,y_1)+d_n(y_1,y_2)+d_n(y_2,x_2) \leq 3$.
\end{center}
Therefore, $\mathrm{diam}_n(Q) \leq 3$.

\medskip \noindent
\underline{Case 2:} Neither $\mathcal{H}$ nor $\mathcal{V}$ contains $n$ pairwise transverse hyperplanes. Naturally, $\mathcal{H}$ (resp. $\mathcal{V}$) can be written as the disjoint union $\mathcal{H}= \mathcal{H}_1 \sqcup \cdots \sqcup \mathcal{H}_r$ (resp. $\mathcal{V}= \mathcal{V}_1 \sqcup \cdots \sqcup \mathcal{V}_s$) where $\mathcal{H}_i$ (resp. $\mathcal{V}_i$) is a set of $\mathrm{Ram}(n)$ consecutive hyperplanes for all $1 \leq i \leq r$ (resp. for all $1 \leq i \leq s$), and $\# \mathcal{H}_r \leq \mathrm{Ram}(n)$ (resp. $\# \mathcal{V}_s \leq \mathrm{Ram}(n)$). Notice that, since $Q$ is $\mathrm{Ram}(n)$-thick, necessarily $r,s \geq 2$ hold. Now, for every $1 \leq i \leq r-1$ and $1 \leq j \leq s-1$, $\mathcal{H}_i$ and $\mathcal{V}_j$ contain $n$ pairwise disjoint hyperplanes; in particular, each of these hyperplanes belongs to an $(n,n)$-grid of hyperplanes, and so is not $n$-combinatorially contracting. Therefore, for every vertex $a \in Q$, there exists a vertex $v(a)$ belonging to a hyperplane $V(a) \in \mathcal{V}$ which is not $n$-combinatorially contracting, such that $d(a,v(a)) \leq 2 \mathrm{Ram}(n)$. Let us fix a hyperplane $H \in \mathcal{H}$ which is not $n$-combinatorially contracting, and for every vertex $a \in Q$, let $w(a)$ be a vertex in $V(a) \cap H$. Finally, for every $x,y \in Q$, we get
\begin{center}
$\begin{array}{lcl} d_n(x,y) & \leq & d_n(x,v(x))+ d_n(v(x),w(x))+d_n(w(x),w(y))+d_n(w(y),v(y))+d_n(v(y),y) \\ \\ & \leq & 2 \mathrm{Ram}(n)+ 1+1+1+2 \mathrm{Ram}(n)  =4 \mathrm{Ram}(n)+3 \end{array}$
\end{center}
Thus, $\mathrm{diam}_n(Q) \leq 4 \mathrm{Ram}(n)+3$. $\square$

\medskip \noindent
Thus, looking at the induced action on the contracting graph, as a consequence of Proposition \ref{contracting} we obtain:

\begin{cor}\label{whyp}
Let $G$ be a group acting geometrically on a CAT(0) cube complex $X$. Then $G$ is weakly hyperbolic relatively to the stabilizers of the non-contracting hyperplanes of $X$.
\end{cor}

\noindent
Recall that a group $G$ is \emph{weakly hyperbolic relatively to a collection of subgroups $\mathcal{H}=\{ H_1, \ldots, H_n \}$} if $G$ acts by isometries on a graph $\Gamma$ such that:
\begin{itemize}
	\item $\Gamma$ is hyperbolic,
	\item $\Gamma$ contains finitely-many orbits of edges,
	\item each vertex-stabilizer is either finite or contains a finite-index subgroup conjugated to some $H_i$,
	\item any $H_i$ stabilizes a vertex.
\end{itemize}

\noindent
\textbf{Proof of Corollary \ref{whyp}.} Let $Y$ be the graph obtained from $X^{(1)}$ by adding a vertex for each non-contracting hyperplanes $J$ and linking by an edge this vertex with any vertex of $N(J)$: this is the usual cone-off we will introduce in Definition \ref{usual}, which is clearly quasi-isometric to the cone-off of Definition \ref{cone-off}. Thus, because $Y$ is quasi-isometric to the contracting graph $\Gamma_{\infty}X$, we deduce from Proposition \ref{contracting} that $Y$ is hyperbolic. Furthermore, the action $G \curvearrowright X$ induces an action $G \curvearrowright Y$ with finitely-many orbits of edges, and the stabilizer of a vertex is either finite if it belongs to $X$ or it is the stabilizer of a non-contracting hyperplane. $\square$

\medskip \noindent
As an application of Corollary \ref{whyp}, we will be able to deduce a weak relative hyperbolicity of right-angled Coxeter groups in Section \ref{RACG}.

\subsection{Strong relative hyperbolicity}

\noindent
After considering the weak relative hyperbolicity, it is natural to focus on the strong relative hyperbolicity. We recall the definition as introduced by Bowditch in \cite{relativelyhyperbolic}.

\begin{definition}
A finitely-generated group $G$ is \emph{(strongly) hyperbolic relatively to a collection of subgroups $\mathcal{H}=\{ H_1, \ldots, H_n \}$} if $G$ acts by isometries on a graph $\Gamma$ such that:
\begin{itemize}
	\item $\Gamma$ is hyperbolic,
	\item $\Gamma$ contains finitely-many orbits of edges,
	\item each vertex-stabilizer is either finite or contains a finite-index subgroup conjugated to some $H_i$,
	\item any $H_i$ stabilizes a vertex,
	\item $\Gamma$ is \emph{fine}, ie., any edge belongs only to finitely-many simple loops (or \emph{cycle}) of a given length.
\end{itemize}
A subgroup conjugated to some $H_i$ is \emph{peripheral}. $G$ is just said \emph{relatively hyperbolic} if it is relatively hyperbolic with respect to a finite collection of proper subgroups.
\end{definition}

\noindent
Thus, a natural question is: when is a cone-off of a CAT(0) cube complex fine? With the definition of a cone-off we used in the previous sections, it is not difficult to notice that essentially a cone-off will never be fine. So we first have to modify our definition:

\begin{definition}\label{usual}
Let $X$ be a CAT(0) cube complex and $\mathcal{Q}$ a collection of subcomplexes. The \emph{usual cone-off of $X$ over $\mathcal{Q}$} is the graph obtained from $X^{(1)}$ by adding a vertex for each subcomplex $Q \in \mathcal{Q}$ and linking it by an edge $Q$ to each vertex belonging to $Q$.
\end{definition}

\noindent
Now, we are able to prove the following criterion. 

\begin{thm}\label{fine}
Let $X$ be a uniformly locally finite CAT(0) cube complex and $Y$ a usual cone-off of $X$ over a collection $\mathcal{Q}$ of combinatorially convex subcomplexes. If $\mathcal{Q}$ is locally finite (ie., there exist only finitely-many subcomplexes of $\mathcal{Q}$ containing a given edge of $X$) and if there exists a constant $C \geq 0$ such that two subcomplexes of $\mathcal{Q}$ are both intersected by at most $C$ hyperplanes, then $Y$ is fine. Conversely, if $\mathcal{Q}$ is not locally finite or if it contains two subcomplexes both intersected by infinitely many hyperplanes, then $Y$ is not fine. 
\end{thm}

\noindent
\textbf{Proof.} Suppose that $\mathcal{Q}$ is locally finite and that there exists a constant $C \geq 0$ such that two subcomplexes of $\mathcal{Q}$ are both intersected by at most $C$ hyperplanes. Let $e \in Y$ be an edge and fix one of its endpoints $a \in X$. To a given cycle $\gamma \subset Y$ of length $n$ containing $e$, we associate a loop $\bar{\gamma} \subset X$ containing $a$ in the following way: The cycle $\gamma$ passes through a sequence of cones $C_1, \ldots, C_k$. For every $1 \leq j \leq k$, let $x_j,y_j \in X$ be the two vertices of $\gamma \cap C_j \cap X$, and choose a combinatorial geodesic $[x_j,y_j]$ between $x_j$ and $y_j$ in $X$. Also, for every $1 \leq j \leq k-1$, choose a combinatorial geodesic $[y_j,x_{j+1}]$ between $y_j$ and $x_{j+1}$ in $X$. Finally, choose two combinatorial geodesics $[a,x_1]$ and $[y_n,a]$ in $X$, respectively between $a$ and $x_1$, and $y_n$ and $a$. Now we set
\begin{center}
$\bar{\gamma}= [a,x_1] \cup [x_1,y_1] \cup \cdots \cup [x_k,y_k] \cup [y_k,a]$.
\end{center}
Fix some $1 \leq j \leq k$. Notice that a hyperplane separating $x_j$ and $y_j$ must intersect $\bar{\gamma} \backslash [x_j,y_j]$, ie., must separate either $a$ and $x_1$, or $a$ and $y_1$, or $x_i$ and $y_i$ for some $1 \leq i \neq j \leq k$, or $y_i$ and $x_{i+1}$ for some $1 \leq i \leq k-1$. But a hyperplane separating $x_j$ and $y_j$ and, for some $1 \leq i \leq k$, $x_i$ and $y_i$, necessarily intersects both $C_j$ and $C_i$. Because there exist at most $C$ such hyperplanes, we deduce that
\begin{center}
$\begin{array}{lcl} d_X(x_j,y_j) & \leq & Ck + d_X(a,x_1)+d_X(a,y_k)+ \sum\limits_{j=1}^{k-1} d_X(y_j,x_{j+1}) \\ \\ & \leq & Cn+n= (C+1)n \end{array}$
\end{center}
Therefore, we get
\begin{center}
$\begin{array}{lcl} \mathrm{diam}_X (\bar{\gamma}) & \leq & d_X(a,x_1) + \sum\limits_{j=1}^{k} d_X(x_j,y_j) + \sum\limits_{j=1}^{k-1} d_X(y_j,x_{j+1}) + d_X(a,x_1) \\ \\ & \leq & k(C+1)n + n \leq (C+2)n^2 \end{array}$ 
\end{center}
We have proved that $\bar{\gamma}$ is included into the ball $B=B_X(a,(C+2)n^2)$. Let $\dot{B}$ denote the cone-off of this ball over $\{Q \cap B \mid Q \in \mathcal{Q} \}$, so that $\bar{\gamma} \subset B$ implies $\gamma \subset \dot{B}$. It is worth noticing that $\dot{B}$ depends only on the edge $e$, the constant $C$ and the integer $n$, so $\dot{B}$ contains all the cycles of length $n$ passing through $e$. Since $\dot{B}$ is finite by the local finiteness of $X$ and $\mathcal{Q}$, we conclude that there exist only finitely many such cycles. This proves that $Y$ is fine.

\medskip \noindent
Conversely, if $\mathcal{Q}$ is not locally finite, then there exist an edge $e \in X$ and infinitely-many subcomplexes $Q_1,Q_2, \ldots \in \mathcal{Q}$ containing $e$. For each $i \geq 1$, $e$ belongs to a cycle of length three passing through the cone associated to $Q_i$. Therefore, $Y$ is not fine. Now, suppose that $\mathcal{Q}$ contains two subcomplexes $C_1,C_2$ both intersected by infinitely-many hyperplanes of $X$. Let $p : X \to C_1$ denote the combinatorial projection onto $C_2$. According to Proposition \ref{proj}, $p(C_2)$ is a geodesic subcomplex of $C_1$ containing infinitely-many hyperplanes, say $J_1,J_2, \ldots$. Choose a basepoint $x \in p(C_2)$ and, for every $i \geq 1$, fix a vertex $x_i \in J \subset p(C_2)$ and a combinatorial geodesic $[x,x_i] \subset p(C_2)$. Because $X$ is locally finite, the sequence of geodesics $([x,x_i])$ contains a subsequence converging to some infinite ray $r : [0,+ \infty) \to p(C_2)$. Because $X$ is finite dimensional, we deduce from the infinite Ramsey Theorem that the set of hyperplanes intersected by $r$ must contain an infinite sequence of pairwise disjoint hyperplanes, say $V_1,V_2,\ldots$; notice that, according to Proposition \ref{proj}, the $V_i$'s intersect $C_1$ as well as $C_2$. For every $n \geq 1$, let $\mathcal{C}_n$ denote the cycle of combinatorially convex subcomplexes $(N(V_1),C_1,N(V_n),C_2)$. From Corollary \ref{flat rectangle}, we get a flat rectangle $D_n$ bounded by $\mathcal{C}$. Once again, because $X$ is locally finite, the sequence of subcomplexes $(D_n)$ contains a subsequence converging to some subcomplexes $D_{\infty}$. Notice that $D_{\infty} \cap C_i$ is an infinite ray $r_i$ (with $i=1,2$) and that $D_{\infty} \cap N(V_i)$ is a combinatorial geodesic $\gamma_i$ between $r_1$ and $r_2$ whose length $\ell$ does not depend on $i$ (with $i \geq 1$). For every $n \geq 2$, let $c_n \subset Y$ denote the cycle passing through $\gamma_1, \gamma_n$ and the two cones over $C_1$ and $C_2$. If $e$ is any edge of $\gamma_1$, then $c_n$ is a sequence of infinitely-many cycles with the same length containing $e$. Therefore, $Y$ is not fine. $\square$

\medskip \noindent
Now, we would like to combine Theorem \ref{KFR} and Theorem \ref{fine} to deduce some criterion of relative hyperbolicity. Although it seems to be difficult to state a general criterion, this sketches a general approach to study the strong relative hyperbolicity of cubulable groups. Roughly speaking:

\medskip \noindent
Let $G$ be a group acting geometrically on a CAT(0) cube complex $X$. If we are able to find a $G$-equivariant collection of combinatorially convex subcomplexes $\mathcal{C}$ such that any thick flat rectangle of $X$ is included into a subcomplex of $\mathcal{C}$, and two subcomplexes of $\mathcal{C}$ do not fellow-travel, then, by looking at the induced action of $G$ on the usual cone-off of $X$ over the collection $\mathcal{Q}$, we should be able to deduce that $G$ is strongly hyperbolic relatively to $\{ \mathrm{stab}(Q) \mid Q \in \mathcal{Q} \}$. In order to construct $\mathcal{Q}$, a possibility could be to start with the collection $\mathcal{F}^0$ of the ``combinatorial flats'' of $X$, and then to define inductively the sequence $(\mathcal{F}^n)$ by: if $C_1, \ldots, C_k$ denote the connected components of the graph whose set of vertices is $\mathcal{F}^n$ and whose edges link two subcomplexes which fellow-travel, we set $\mathcal{F}^{n+1}= \left( \mathrm{co} \left( \bigcup\limits_{ F \in C_1} F \right), \ldots , \mathrm{co} \left( \bigcup\limits_{F \in C_k} F \right) \right)$, where $\mathrm{co}(\cdot)$ denotes the combinatorial convex hull. Finally, because $X$ is cocompact, the sequence $(\mathcal{F}^n)$ should be eventually constant to some collection $\mathcal{F}^{\infty}$, and this is our condidate for $\mathcal{Q}$. 

\medskip \noindent
A similar idea can be found in \cite{isolated}, and is successfully applied to Coxeter groups \cite[Theorem A.1]{BHSC}. In the next section, following the argument we have sketched above, we will be able to give a direct proof of the second part of \cite[Theorem I]{BHSC}, characterizing the strong relative hyperbolicity of right-angled Coxeter groups.

\subsection{Application to right-angled Coxeter groups}\label{RACG}

\noindent
Our plan is to apply the results established in the two previous sections to the class of right-angled Coxeter groups. We begin with some basic definitions.

\begin{definition}
Let $\Gamma$ be a simplicial graph. The \emph{right-angled Coxeter group} $C(\Gamma)$ is defined by the presentation
\begin{center}
$\langle v \in V(\Gamma) \mid u^2=1, \ [v,w]=1, \ u \in V(\Gamma), (v,w) \in E(\Gamma) \rangle$,
\end{center}
where $V(\Gamma)$ and $E(\Gamma)$ denote respectively the sets of vertices and edges of $\Gamma$.
\end{definition}

\noindent
If $\Lambda \leq \Gamma$ is an \emph{induced subgraph} (ie., two vertices of $\Lambda$ are adjacent in $\Lambda$ if and only if they are adjacent in $\Gamma$), then the subgroup $\langle \Lambda \rangle$ generated by the vertices of $\Lambda$ is naturally isomorphic to the right-angled Coxeter group $C(\Lambda)$. For instance, if $\Lambda$ is a \emph{clique} of $\Gamma$ (ie., a complete subgraph) with $n$ vertices, then $\langle \Lambda \rangle$ is isomorphic to direct product of $n$ copies of $\mathbb{Z}_2$: such a subgroup is referred to as a \emph{clique subgroup}. Another interesting class of subgroups is:

\begin{definition}
The \emph{join} $\Lambda_1 \ast \Lambda_2$ of two graphs $\Lambda_1$ and $\Lambda_2$ is the graph obtained from the disjoint union $\Lambda_1 \sqcup \Lambda_2$ by adding an edge between any vertex of $\Lambda_1$ and any vertex of $\Lambda_2$. By extension, if a given simplicial graph $\Gamma$ contains a subgraph $\Lambda$ which splits as a join of two non empty subgraphs, we say that the subgroup $\langle \Lambda \rangle$ of the right-angled Coxeter group $C(\Gamma)$ is a \emph{join subgroup}.
\end{definition}

\noindent
Given a simplicial graph $\Gamma$, the right-angled Coxeter group $C(\Gamma)$ naturally acts on a CAT(0) cube complex $X(\Gamma)$ whose:
\begin{itemize}
	\item vertices are the elements of $C(\Gamma)$,
	\item edges link two vertices $g,h$ if $h=gv$ for some $v \in V(\Gamma)$,
	\item $n$-cubes are generated by the vertices 
\begin{center}
$\{ gv_{i_1}\cdots v_{i_r} \mid 0 \leq r \leq n, 1 \leq i_1 < \cdots < i_r \leq n \}$,
\end{center}
				where $g \in C(\Gamma)$ and $v_1, \ldots, v_n \in V(\Gamma)$ are pairwise adjacent.
\end{itemize}
Notice that the 1-skeleton of $X(\Gamma)$ is the Cayley graph of $C(\Gamma)$ with respect to the generating set $V(\Gamma)$. In particular, if $J_v$ denotes the hyperplane dual to the edge $(1,v)$ where $v \in V(\Gamma)$, then the hyperplanes of $X(\Gamma)$ are precisely the translates of the $J_v$'s. 

\medskip \noindent
The two following lemmas on the geometry of $X(\Gamma)$ are elementary, and their proofs are left to the reader. In the case of right-angled Artin groups, the analogous results have been proved in \cite[Section 3]{BehrstockCharney}.

\begin{lemma}\label{transverse}
If the hyperplanes $gJ_a$ and $hJ_b$ are transverse, then $a \neq b$ and $[a,b]=1$.
\end{lemma}

\begin{lemma}\label{stab}
The stabilizer of the hyperplane $gJ_a$ is $g \langle \mathrm{star}(a) \rangle g^{-1}$. 
\end{lemma}

\noindent
Recall that, given a graph $\Gamma$ and a vertex $v$, the \emph{link} of $v$, denoted by $\mathrm{link}(v)$, is defined as the subgraph generated by the vertices of $\Gamma$ adjacent to $v$; and its \emph{star}, denoted by $\mathrm{star}(v)$, is defined as the subgraph generated by $\{v\} \cup \mathrm{link}(v)$.

\medskip \noindent
We begin our study of the relative hyperbolicity of right-angled Coxeter groups by determining precisely when a hyperplane of our cube complex is contracting.

\begin{prop}\label{hypcon}
Let $\Gamma$ be a finite graph and $u \in V(\Gamma)$. The hyperplane $J_u$ of $X(\Gamma)$ is contracting if and only if $u$ does not belong to an induced square in $\Gamma$. 
\end{prop}

\noindent
\textbf{Proof.} Suppose that the vertices $u,v,w,x \in V(\Gamma)$ generate a square in $\Gamma$ (with $u$ and $x$ both adjacent to $v$ and $w$). Then $X(\Gamma)$ contains the following subcomplex for every $n \geq 1$:
\begin{center}
\includegraphics[scale=0.6]{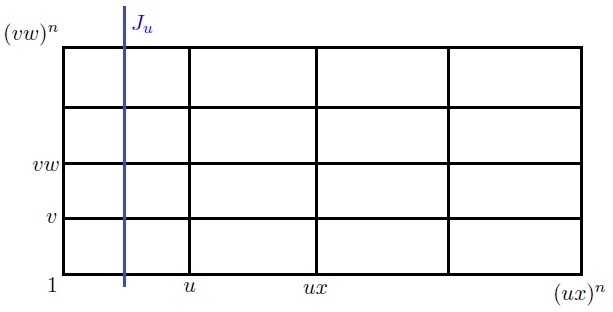}
\end{center}
It follows from Lemma \ref{transverse} that the hyperplanes 
\begin{center}
$(\{J_u,uJ_x,uxJ_u,\ldots, (ux)^{n-1}uJ_x\}, \{ J_v, vJ_w, vwJ_v, \ldots, (vw)^{n-1}vJ_w\})$
\end{center}
define a $(2n,2n)$-grid of hyperplanes. Thus, $J_u$ is not contracting.

\medskip \noindent
Conversely, suppose that $J_u$ is not contracting. In particular, $J_u$ belongs to a $(L,L)$-grid of hyperplanes $(\mathcal{H}, \mathcal{V})$ satisfying $L >n$, where $n$ is the maximal cardinality of a clique subgroup of $C(\Gamma)$. To fix the notation, say that $\mathcal{H}= \{ A_1, \ldots, A_r, J_u, B_1, \ldots, B_s\}$ and $\mathcal{V} = \{ V_1, \ldots, V_{L}\}$, with the convention that each of these hyperplanes separates its two adjacent hyperplanes in the list it belongs. Without loss of generality, we may suppose that $s \geq 1$. Let $J_u^-$ (resp. $B_s^+$) denote the half-space delimited by $J_u$ (resp. $B_s$) not containing $B_s$ (resp. $J_u$) and $V_1^-$ (resp. $V_{L}^+$) the half-space delimited by $V_1$ (resp. $V_{L}$) not containing $V_{L}$ (resp. $V_1$). According to Corollary \ref{flat rectangle}, a disc diagram of minimal complexity bounded by the cycle of subcomplexes $(J_u^-,V_{L}^+, B_s^+,V_1^-)$ will define a subcomplex $Q$ isomorphic to a rectangle. By translating $Q$ by an element of $\langle \mathrm{star}(u) \rangle$ if necessary, we may suppose without loss of generality that $1$ is the corner of $Q$ in $V_1^- \cap J_u^-$. Let $p_1 \cdots p_{i}$ denote the corner of $Q$ in $J_u^- \cap V_{2L}^+$ and $uq_1 \cdots q_{j}$ the corner of $Q$ in $V_1^- \cap B_s^+$, where $p_1, \ldots, p_{i},q_1, \ldots, q_{j} \in V(\Gamma)$. 

\medskip \noindent
Suppose by contradiction that $u$ is adjacent to $q_1, \ldots, q_j$ in $\Gamma$. Let $B_1$ be the hyperplane dual to the edge $(uq_1 \cdots q_k, uq_1 \cdots q_{k+1})$, ie., $B_1= uq_1 \cdots q_k J_{q_{k+1}}$. Now, notice that $J_u$ and $B_1$ are transverse if and only if $q_k \cdots q_1uJ_u=J_u$ and $J_{q_{k+1}}$ are transverse if and only if $u$ and $q_{k+1}$ are adjacent vertices in $\Gamma$. On the other hand, $J_u$ and $B_1$ are disjoint by assumption, so there exists $1 \leq h \leq s$ such that $u$ and $q_h$ are not adjacent in $\Gamma$.

\medskip \noindent
Moreover, we deduce from the definition of $n$ that $\langle p_1, \ldots, p_L \rangle$ is an infinite subgroup, so that there exist $1 \leq \ell < m \leq L$ such that $p_{\ell}$ and $p_m$ are not adjacent in $\Gamma$. 

\medskip \noindent
Finally, noticing that $u,q_1, \ldots, q_j$ are all adjacent to $p_1, \ldots, p_i$ according to Lemma \ref{transverse}, we conclude that $u$ belongs to the induced square generated by $u,q_h,p_{\ell},p_m$ in $\Gamma$. $\square$

\medskip \noindent
Thus, if $\square(\Gamma)$ denote the set of the vertices of a graph $\Gamma$ which belong to an induced square, we obtain by combining Corollary \ref{whyp}, Proposition \ref{hypcon} and Lemma \ref{stab}:

\begin{prop}\label{RACG1}
Let $\Gamma$ be a finite graph. The right-angled Coxeter group $C(\Gamma)$ is weakly hyperbolic relatively to $\{ \langle \mathrm{star}(u) \rangle \mid u \in \square(\Gamma) \}$.
\end{prop}

\noindent
In fact, during the proof of Proposition \ref{hypcon}, we have shown:

\begin{fact}\label{flatandjoin}
Let $\Gamma$ be a finite graph. If $n$ denotes the maximal cardinality of a clique subgroup of $C(\Gamma)$, then any $n$-thick flat rectangle in $X(\Gamma)$ is included into a subcomplex $X(\Gamma_1 \ast \Gamma_2) \subset X(\Gamma)$, where $\Gamma_1,\Gamma_2 \subset \Gamma$ are not complete subgraphs. 
\end{fact}

\noindent
Therefore, combined with Theorem \ref{KFR}, we deduce that right-angled Coxeter groups are also weakly hyperbolic relatively to their join subgroups.

\begin{prop}\label{RACG2}
Let $\Gamma$ be a finite graph. The right-angled Coxeter group $C(\Gamma)$ is weakly hyperbolic relatively to $\{ C(\Gamma_1 \ast \Gamma_2) \mid \Gamma_1 \ast \Gamma_2 \subset \Gamma \ \text{where} \ \Gamma_1,\Gamma_2 \ \text{are not complete}\}$. 
\end{prop}

\begin{remark}
Proposition \ref{RACG1} and Proposition \ref{RACG2} should be compared with the weak relative hyperbolicity established in \cite{wrelhypRACG}.
\end{remark}

\noindent
From now on, let us focus on the strong relative hyperbolicity. We begin with some preliminary definitions.

\begin{definition}
Let $\Gamma$ be a simplicial graph. A join $\Gamma_1 \ast \Gamma_2$ in $\Gamma$ is \emph{large} if neither $\Gamma_1$ nor $\Gamma_2$ is a complete subgraph. 
\end{definition}

\begin{definition}
Let $\Gamma$ be a finite graph. A \emph{join decomposition} of $\Gamma$ is a collection of subgraphs $(\Gamma_1, \ldots, \Gamma_n)$, with $n \geq 0$, such that:
\begin{itemize}
	\item any large join of $\Gamma$ is included into $\Gamma_i$ for some $1 \leq i \leq n$,
	\item $\Gamma_i \cap \Gamma_j$ is complete for every $1 \leq i < j \leq n$,
	\item for every vertex $v$, $\mathrm{link}(v) \cap \Gamma_i$ not complete implies $v \in \Gamma_i$.
\end{itemize}
There exists at least one join decomposition: the \emph{trivial} decomposition $(\Gamma)$. 
\end{definition}

\noindent
According to the following proposition, join decompositions lead to strong relative hyperbolicity.

\begin{prop}\label{rel}
Let $\Gamma$ be a finite graph and $(\Gamma_1, \ldots, \Gamma_n)$ a join decomposition. Then the right-angled Coxeter group $C(\Gamma)$ is hyperbolic relatively to $\{C(\Gamma_1), \ldots, C(\Gamma_n) \}$.
\end{prop}

\noindent
We first need to prove the following lemma:

\begin{lemma}\label{sub}
Let $X$ be a uniformly locally finite CAT(0) cube complex. For every $N \geq 1$, there exists a constant $M \geq 1$ such that, for any collection $\mathcal{H}$ of at least $M$ hyperplanes, $\mathcal{H}$ contains two hyperplanes separated by at least $N$ hyperplanes of $X$. 
\end{lemma}

\noindent
\textbf{Proof.} Let $J_1, \ldots, J_r$ be $r$ hyperplanes such that, for every $1 \leq i \neq j \leq r$, $J_i$ and $J_j$ are not separated by $N$ hyperplanes. For every $1 \leq i \leq r$, let $N_i$ denote the $N$-neighborhood of $J_i$ with respect to $d_{\infty}$ (see Section \ref{d-infty}); according to \cite[Corollary 3.5]{packing}, $N_i$ is combinatorially convex. Furthermore, our assumption implies that $N_i \cap N_j \neq \emptyset$ for every $1 \leq i , j \leq r$, hence $\bigcap\limits_{i=1}^r N_i \neq \emptyset$. Noticing that
\begin{center}
$d_{\infty} \leq d \leq \dim(X) \cdot d_{\infty}$,
\end{center}
we deduce that, if we fix some vertex $c \in \bigcap\limits_{i=1}^r N_i$, the hyperplanes $J_1, \ldots, J_r$ intersect the ball $B(c,N \cdot \dim(X))$. On the other hand, $X$ is uniformly locally finite, so the cardinality of this ball is bounded above by a constant depending only on $N$, and a fortiori the number of hyperplanes intersecting this ball is bounded above by a constant $\kappa(N)$. Therefore, setting $M= \kappa(N)+1$ proves our lemma. $\square$

\medskip \noindent
\textbf{Proof of Proposition \ref{rel}.} Let $Y$ denote the cone-off of $X(\Gamma)$ over the translates of the $X(\Gamma_i)$'s. Combining Fact \ref{flatandjoin} and Theorem \ref{KFR}, we know that $Y$ is hyperbolic. To conclude, it is sufficient to deduce from Theorem \ref{fine} that $Y$ is fine. 

\medskip \noindent
First, we will show that this collection of subcomplexes is locally finite. Indeed, if an edge $(g,gu)$ belongs to pairwise distinct subcomplexes $h_1C(\Lambda_1), h_2C(\Lambda_2), \ldots$, where $\Lambda_i \in \{\Gamma_1, \ldots, \Gamma_n \}$ for $i \geq 1$, then $u \in \bigcap\limits_{i \geq 1} \Lambda_i$ and $g \in \bigcap\limits_{i \geq 1} h_iC(\Lambda_i)$. Since two different cosets of a given subgroup are necessarily disjoint, we deduce that the $\Lambda_i$'s have to be pairwise distinct, otherwise the intersection $\bigcap\limits_{i \geq 1} h_i C(\Lambda_i)$ would be empty. Because we have only $n$ subgraphs in our join decomposition, this means that our collection $h_1C(\Lambda_1), h_2C(\Lambda_2), \ldots$ is necessarily finite. 

\medskip \noindent
Now, we want to prove the separation property. Let $N$ denote the maximal order of a finite subgroup of $C(\Gamma)$ and let $M$ be the constant given by Lemma \ref{sub}. Suppose that there exist $M$ hyperplanes intersecting both $gX(\Gamma_i)$ and $h X(\Gamma_j)$. By definition of $M$, there exist two hyperplanes $V_1,V_2$, intersecting both $gX(\Gamma_i)$ and $h X(\Gamma_j)$, separated by $N$ hyperplanes. Therefore, since $\mathcal{C}=(V_1,gX(\Gamma_i),V_2,hX(\Gamma_j))$ is a cycle of four combinatorially convex subcomplexes, we deduce from Corollary \ref{flat rectangle} that there exists a flat rectangle $D$ bounded by $\mathcal{C}$. Let $g_1,g_2,g_3,g_4$ denote the corners of $D$ cyclically ordered, with $g_1,g_2 \in h C(\Gamma_j)$ and $g_3,g_4 \in g C(\Gamma_i)$. Set $w_1=g_1^{-1}g_2$ and $w_2=g_1^{-1}g_4$, so that $g_2=g_1w_1$ and $g_4=g_1w_2$; notice that $w_1$ and $w_2$ commute, the product $w_1w_2=w_2w_1$ is reduced in $C(\Gamma)$ and $g_3=g_1w_1w_2=g_1w_2w_1$. Let $\mathrm{supp}(w_i)$ denote the set of letters (or equivalently, vertices) used to write the word $w_i$ ($i=1,2$). Because $V_1$ and $V_2$ are separated by $N$ hyperplanes, we deduce that $\langle \mathrm{supp}(w_1) \rangle$ has cardinality at least $N$; by the definition of $N$, this implies that $\langle \mathrm{supp}(w_1) \rangle$ is infinite, ie., $\mathrm{supp}(w_1)$ contains two non-adjacent vertices $v_1,v_2$. 

\medskip \noindent
From $g_1,g_2 \in h C(\Gamma_j)$, we deduce that $w_1=g_1^{-1}g_2 \in C(\Gamma_j)$; from $g_3,g_4 \in g C(\Gamma_i)$, we deduce that $w_1=g_4^{-1}g_3 \in C(\Gamma_i)$. Therefore, 
\begin{center}
$\{v_1,v_2 \} \subset \mathrm{supp}(w_1) \subset \Gamma_i \cap \Gamma_j$. 
\end{center}
This implies $i=j$. 

\medskip \noindent
Because $w_2$ commute with $w_1$ and because the product $w_1w_2$ is reduced in $C(\Gamma)$, any vertex $v \in \mathrm{supp}(w_2)$ satisfies 
\begin{center}
$\mathrm{link}(v) \cap \Gamma_j \supset \mathrm{link}(v) \cap \mathrm{supp}(w_1) \supset \{ v_1,v_2 \}$.
\end{center}
This implies $v \in \Gamma_j$. A fortiori, we get $w_2 \in \Gamma_j$ hence $g_4 =g_1w_2 \in hC(\Gamma_j)=hC(\Gamma_i)$. Thus,
$g_4 \in gC(\Gamma_i) \cap hC(\Gamma_i)$, so $gC(\Gamma_i)=hC(\Gamma_i)$. 

\medskip \noindent
We have proved that if two subcomplexes of our collection are both intersected by $N$ hyperplanes then they are equal. This concludes the proof. $\square$

\medskip \noindent
Although finding a non trivial join decomposition implies that our right-angled Coxeter group is relatively hyperbolic, a priori we do not know how to find such a decomposition, and if so some decompositions are better than others. For instance, if a graph $\Gamma$ contains cut-vertices, then any choice of a set $S$ of these cut-vertices produces a join decomposition $(\Gamma_1, \ldots, \Gamma_r)$, where, for every $i \neq j$, $\Gamma_i \cap \Gamma_j$ is either empty or a single vertex of $S$; clearly, the decomposition seems to be ``finer'' when $S$ is the set of all the cut-points of $\Gamma$. This is the motivation to introduce the following canonical join decomposition of a graph:

\begin{definition}\label{decomposition}
Let $\Gamma$ be a finite graph. For every subgraph $\Lambda \subset \Gamma$, let $\mathrm{cp}(\Lambda)$ denote the subgraph of $\Gamma$ generated by $\Lambda$ and the vertices $v \in \Gamma$ such that $\mathrm{link}(v) \cap \Lambda$ is not complete. Now, define the collection of subgraphs $\mathfrak{J}^n(\Gamma)$ of $\Gamma$ by induction in the following way:
\begin{itemize}
	\item $\mathfrak{J}^0(\Gamma)$ is the collection of all the large joins in $\Gamma$;
	\item if $C_1, \ldots, C_k$ denote the connected components of the graph whose set of vertices is $\mathfrak{J}^n(\Gamma)$ and whose edges link two subgraphs with non-complete intersection, we set $\mathfrak{J}^{n+1}(\Gamma)= \left( \mathrm{cp} \left( \bigcup\limits_{ \Lambda \in C_1} \Lambda \right), \ldots , \mathrm{cp} \left( \bigcup\limits_{\Lambda \in C_k} \Lambda \right) \right)$.
\end{itemize}
Because $\Gamma$ is finite, the sequence $(\mathfrak{J}^n(\Gamma))$ must be eventually constant to some collection $\mathfrak{J}^{\infty}(\Gamma)$; see Example \ref{exjoindecomposition}. By construction, it is clear that $\mathfrak{J}^{\infty}(\Gamma)$ is a join decomposition of $\Gamma$. Furthermore, according to the next lemma, it may be thought of as the minimal join decomposition. 
\end{definition}

\begin{lemma}
Let $\Gamma$ be a finite graph and $(\Gamma_1, \ldots, \Lambda_m)$ a join decomposition. For every $\Lambda \in \mathfrak{J}^{\infty} (\Gamma)$, there exists some $1 \leq i \leq m$ such that $\Lambda \subset \Gamma_i$. 
\end{lemma}

\noindent
\textbf{Proof.} We prove the statement by induction on $n$. If $n=0$, this is true by the definition of a join decomposition. Now, suppose that the statement is true for $\mathfrak{J}^n(\Gamma)$ and let $\Lambda \in \mathfrak{J}^{n+1}(\Gamma)$. The subgraph $\Lambda$ corresponds to a connected component $C=\{\Lambda_1, \ldots, \Lambda_k\}$ of the graph whose set of vertices is $\mathfrak{J}^n(\Gamma)$ and whose edges link two subgraphs with non-complete intersection, ie., $\Lambda= \mathrm{cp}(\Lambda_1 \cup \cdots \cup \Lambda_k)$. By our induction hypothesis, for every $1 \leq s \leq k$, there exists $1 \leq i_s \leq m$ such that $\Lambda_s \subset \Gamma_{i_s}$. Noticing that, for every $1 \leq r < s \leq k$, $\Gamma_{i_k} \cap \Gamma_{i_s}$ contains the non-complete subgraph $\Lambda_{i_k} \cap \Lambda_{i_s}$, we deduce that $i_k=i_s$; let $\iota$ be this common value. We have
\begin{center}
$\Lambda = \mathrm{cp}(\Lambda_1 \cup \cdots \cup \Lambda_k) \subset \mathrm{cp}(\Gamma_{\iota})=\Gamma_{\iota}$.
\end{center}
This concludes the proof. $\square$

\begin{remark}
In Definition \ref{decomposition}, setting $\mathfrak{J}^0(\Gamma)$ as the collection of the induced squares in $\Gamma$ does not affect $\mathfrak{J}^{\infty}(\Gamma)$. In practice, it may be simpler to determine the decomposition $\mathfrak{J}^{\infty}(\Gamma)$ in this way. Moreover, it turns out that the subgraphs of $\mathfrak{J}^{\infty}(\Gamma)$ are precisely the maximal subgraphs of $\Gamma$ which belong to the class $\mathcal{T}$ introduced in \cite{BHSC}.
\end{remark}

\noindent
Finally, we are able to characterize the strong relative hyperbolicity of right-angled Coxeter groups, reproving the second part of \cite[Theorem I]{BHSC}.

\begin{thm}\label{mainthrelathypRACG}
Let $\Gamma$ be a finite graph. The right-angled Coxeter group $C(\Gamma)$ is relatively hyperbolic if and only if $\mathfrak{J}^{\infty} (\Gamma)$ is not the trivial decomposition. If so, then $C(\Gamma)$ is hyperbolic relatively to the collection $\{ C(\Lambda) \mid \Lambda \in \mathfrak{J}^{\infty} (\Gamma) \}$.
\end{thm}

\noindent
\textbf{Proof.} Suppose that $\mathfrak{J}^{\infty}(\Gamma)$ is the trivial decomposition and that $C(\Gamma)$ is hyperbolic relatively to some collection of subgroups $\mathcal{H}$. We will argue by induction on $n$ that, for any $n \geq 0$ and $\Lambda \in \mathfrak{J}^n(\Gamma)$, the subgroup $C(\Lambda)$ is included into a peripheral subgroup. Because any subgroup isomorphic to a direct product of two infinite groups has to be included into a peripheral subgroup, the statement holds for $n=0$. Now suppose that this statement holds for some $n \geq 0$, and let $\Lambda \in \mathfrak{J}^{n+1}(\Gamma)$. The subgraph $\Lambda$ corresponds to a connected component $C=\{\Lambda_1, \ldots, \Lambda_k\}$ of the graph whose set of vertices is $\mathfrak{J}^n(\Gamma)$ and whose edges link two subgraphs with non-complete intersection, ie., $\Lambda= \mathrm{cp}(\Lambda_1 \cup \cdots \cup \Lambda_k)$. By our induction hypothesis, for every $1 \leq i \leq k$, the subgroup $C(\Lambda_i)$ is included into some peripheral subgroup $H_i$. Notice that, for every $1 \leq i <j \leq k$, $H_i \cap H_j$ contains $C(\Lambda_i) \cap C(\Lambda_j)=C(\Lambda_i \cap \Lambda_j)$, which is infinite, hence $H_i=H_j$. Therefore, the $C(\Lambda_i)$'s are all included into the same peripheral subgroup $H$. Then, for every vertex $v \in \mathrm{cp}(\Lambda_1 \cup \cdots \cup \Lambda_k)$, we have
\begin{center}
$H \cap H^v \supset C(\Lambda_1 \cup \cdots \cup \Lambda_k) \cap C(\Lambda_1 \cup \cdots \cup \Lambda_k)^v =C(\Lambda_1 \cup \cdots \cup \Lambda_k)$,
\end{center}
so $H \cap H^v$ is infinite, and this implies $v \in H$. Therefore,
\begin{center}
$C(\Lambda) = C(\mathrm{cp}(\Lambda_1 \cup \cdots \cup \Lambda_k)) \subset H$.
\end{center}
This concludes the proof of our claim. Finally, because $\mathfrak{J}^{\infty}(\Gamma)=(C(\Gamma))$, we deduce that some peripheral subgroup of $\mathcal{H}$ is not proper. Consequently, $C(\Gamma)$ is not relatively hyperbolic.

\medskip \noindent
The converse is a consequence of Proposition \ref{rel}. $\square$

\begin{remark}
It is worth noticing that the previous proof shows that the collection of subgroups $\{ C(\Lambda) \mid \Lambda \in \mathfrak{J}^{\infty} (\Gamma) \}$ is a minimal collection of peripheral subgroups, ie., if $C(\Gamma)$ hyperbolic relatively to $\mathcal{H}$, then, for every $\Lambda \in \mathfrak{J}^{\infty}(\Gamma)$, there exists $H \in \mathcal{H}$ such that $C(\Lambda) \subset H$.
\end{remark}

\begin{ex}\label{exjoindecomposition}
The figure below gives an example of a graph $\Gamma$ and its canonical join decomposition $\mathfrak{J}^{\infty}(\Gamma)$.
\begin{center}
\includegraphics[scale=0.6]{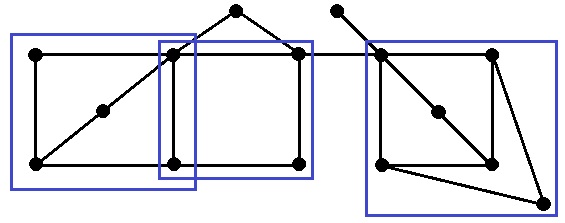}
\end{center}
In particular, we deduce from Theorem \ref{mainthrelathypRACG} that the associated right-angled Coxeter group $C(\Gamma)$ is relatively hyperbolic.
\end{ex}

\section{Hyperbolicity with respect to the \texorpdfstring{$\ell_{\infty}$}{l-infinity} metric}\label{d-infty}

\noindent
Recall that the distance $d_{\infty}$ on a CAT(0) cube complex $X$ is the extension on the whole cube complex of the $\ell_{\infty}$-norm defined on each cube. Alternatively, the restriction of $d_{\infty}$ over $X^{(0)}$ corresponds to the graph metric associated to the cone-off of $X$ over its cubes. Because, loosely speaking, $d_{\infty}$ ``kills'' the dimension of $X$, it may be expected, following Theorem \ref{critère d'hyperbolicité}, that $(X,d_{\infty})$ is hyperbolic precisely when its grids of hyperplanes cannot be too large. This is precisely what we prove.

\begin{thm}\label{dinfty}
Let $X$ be a CAT(0) cube complex. Then $(X,d_{\infty})$ is hyperbolic if and only if the grids of hyperplanes in $X$ are uniformly thin.
\end{thm}

\noindent
We begin by proving a preliminary lemma.

\begin{lemma}\label{lemme}
Let $X$ be a CAT(0) cube complex and $x,y,z,z' \in X$ four pairwise distinct vertices. Let $m=m(x,y,z)$ and $m'=m(x,y,z')$ be the associated median vertices. Then any hyperplane separating $m$ and $m'$ separates $z$ and $z'$ as well. In particular, $d_{\infty}(m,m') \leq d_{\infty}(z,z')$.
\end{lemma}

\noindent
\textbf{Proof.} Fix eight combinatorial geodesics $[x,m]$, $[y,m]$, $[z,m]$, $[x,m']$, $[y,m']$, $[z,m']$, $[m,m']$, $[z,z']$. We get a bigon $B= \{ [x,m] \cup [m,y], [x,m'] \cup [m',y] \}$ and a $4$-gon $P= \{[m,m'],[m',z'],[z',z],[z,m]\}$. Given a hyperplane $J$ separating $m$ and $m'$, in $B$ necessarily either $J$ intersects $[x,m']$ and $[y,m]$ or it intersects $[x,m]$ and $[y,m']$. Say we are in the first case, the second one being completely similar. In particular, $J$ cannot intersect neither $[m,z]$ nor $[m',z']$ since otherwise $J$ would intersect the combinatorial geodesic $[y,m] \cup [m,z]$ or $[x,m'] \cup [m',z']$ twice. On the other hand, $J$ has to intersect the $4$-gon $P$ since it separates $m$ and $m'$, so we conclude that $J$ must separate $z$ and $z'$. $\square$

\medskip \noindent
\textbf{Proof of Theorem \ref{dinfty}.} Suppose that the grids of hyperplanes in $X$ are all $C$-thin for some $C \geq 1$. According to Lemma \ref{thinbigon}, it is sufficient to prove that the bigons in $X$ are all $(C+3)$-thin with respect to $d_{ \infty}$. 

\medskip \noindent
Let $\{ \gamma_1, \gamma_2 \} \subset X$ be a bigon with endpoints $x,y$; for convenience, let $\ell = d_X(x,y)$. Let us fix some vertex $z \in \gamma_1$, and, for any $0 \leq t \leq \ell$, let $z(t)$ denote the vertex of $\gamma_2$ defined by $d_X(x,z(t))=t$. We may suppose that $d_{\infty}(x,z),d_{\infty}(z,y) \geq C+2$, otherwise there is nothing to prove. Finally, let $m_1(t)$ and $m_2(t)$ be the median points associated to $\{x,z,z(t)\}$ and $\{y,z,z(t)\}$ respectively, and let $F$ denote the function $t \mapsto d_{\infty}(z,m_1(t))-d_{\infty}(z,m_2(t))$.
\begin{center}
\includegraphics[scale=0.6]{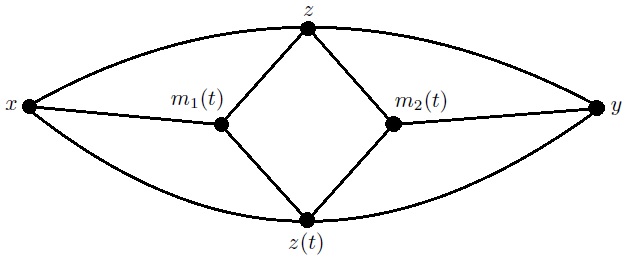} 
\end{center}
Notice that, using Lemma \ref{lemme}, we get
\begin{center}
$d_{\infty}(m_2(1),z=m(x,y,z)) \leq d_{\infty}(x,z(1))=1$.
\end{center}
Similarly, $d_{\infty}(m_1(\ell-1),z) \leq 1$. Futhermore, we have 
\begin{center}
$d_{\infty}(z,m_1(1)) \geq d_{\infty}(z,x)- d_{\infty}(x,m_1(1)) \geq C+2-1=C+1$.
\end{center}
Similarly, $d_{\infty}(z,m_2(\ell-1)) \geq C+1$. Therefore, we conclude that $F(1) \geq C$ and $F(\ell-1) \leq -C$. On the other hand, for all $1 \leq t \leq d(x,y)-1$, we have
\begin{center}
$\begin{array}{lcl} |F(t+1)-F(t)| & =& |d_{\infty}(z,m_1(t+1))-d_{\infty}(z,m_2(t+1))- d_{\infty}(z,m_1(t))+d_{\infty}(z,m_2(t))| \\ \\ & \leq & |d_{\infty}(z,m_1(t+1))-d_{\infty}(z,m_1(t))|+ |d_{\infty}(z,m_2(t+1))-d_{\infty}(z,m_2(t))| \\ \\ & \leq & d_{\infty}(m_1(t),m_1(t+1))+ d_{\infty}(m_2(t),m_2(t+1)) \\ \\ & \leq & d_{\infty}(z(t),z(t+1))+d_{\infty}(z(t),z(t+1))=2 \end{array}$
\end{center}
where the last inequality is justified by Lemma \ref{lemme}. Consequently, there exists some $1 \leq t_0 \leq d(x,y)$ such that $|F(t_0)| \leq 3$. That is to say, if $V_1, \ldots, V_r$ (resp. $H_1, \ldots, H_s$) is a maximal collection of pairwise disjoint hyperplanes separating $z$ and $m_1(t_0)$ (resp. $z$ and $m_2(t_0)$), ie., with $r=d_{\infty}(z,m_1(t_0))$ (resp. $s=d_{\infty}(z,m_2(t_0))$), then $|r-s| \leq 3$. 

\medskip \noindent
Applying Fact \ref{median}, we notice that $z$ is the median point $m(z,m_1(t_0),m_2(t_0))$, and similarly $z(t_0)$ is the median point $m(z(t_0),m_1(t_0),m_2(t_0))$. Thus, no hyperplane can separate $z$ or $z(t_0)$ from $\{m_1(t_0),m_2(t_0)\}$; we deduce that the $V_i$'s separate $\{ z,m_2(t_0)\}$ and $\{ m_1(t_0),z(t_0)\}$, and similarly the $H_j$'s separate $\{ z,m_1(t_0) \}$ and $\{m_2(t_0),z(t_0)\}$. It follows that $V_i$ and $H_j$ are transverse for any $1 \leq i \leq r$ and $1 \leq j \leq s$, and our hypothesis implies that $\min(r,s) \leq 3$. Finally, we conclude that
\begin{center}
$d_{\infty}(z,\gamma_2) \leq d_{\infty}(z,z(t_0)) = \max(r,s) \leq C+3$,
\end{center}
ie., the bigon $\{ \gamma_1,\gamma_2 \}$ is $(C+3)$-thin. 

\medskip \noindent
Conversely, suppose that $(X,d_{\infty})$ is $\delta$-hyperbolic. Let $(\mathcal{H},\mathcal{V})$ be a grid of hyperplanes in $X$. We want to prove that $\min(\# \mathcal{H},\# \mathcal{V}) \leq 4\delta +2$. Of course, if $\# \mathcal{V} \leq 4\delta +3$, there is nothing to prove, so we may suppose without loss of generality that $\# \mathcal{V} > 4\delta +3$. For convenience, let $\mathcal{H}=\{H_1, \ldots, H_r\}$ (resp. $\mathcal{V}= \{ V_1, \ldots, V_s \}$) such that $H_i$ separates $H_{i-1}$ and $H_{i+1}$ for all $1 \leq i \leq r$ (resp. $V_j$ separates $V_{j-1}$ and $V_{j+1}$ for all $1 \leq j \leq s$); let $H_1^-$ and $H_r^+$ (resp. $V_1^-$ and $V_s^+$) be the two disjoint half-spaces delimited respectively by $H_1$ and $H_r$ (resp. $V_1$ and $V_s$). 

\medskip \noindent
Let $x \in H_r^+ \cap V_1^-$, $y \in V_s^+ \cap H_r^+$, $z \in V_s^+ \cap H_1^-$ and $w \in V_1^- \cap H_1^-$ be four vertices, and fix five $d_{\infty}$-geodesics $\gamma(x,y)$, $\gamma(y,z)$, $\gamma(z,w)$, $\gamma(w,x)$ and $\gamma(y,w)$. According to Proposition \ref{hsgeodesic}, we may suppose without loss of generality that $\gamma(x,y) \subset H_r^+$, $\gamma(y,z) \subset V_s^+$, $\gamma(z,w) \subset H_1^-$ and $\gamma(w,x) \subset V_1^-$. Let $p \in \gamma(x,y)$ be a vertex satisfying $|d_{\infty}(x,p)-d_{\infty}(x,y)/2| \leq 1/2$. Notice that
\begin{center}
$d_{\infty}(p,\gamma(x,w)) \geq \frac{s}{2}-2 > \delta$,
\end{center}
since $V_1, \ldots, V_{\lfloor s/2 \rfloor-1}$ separates $p$ and $\gamma(x,w)$. Thus, because the geodesic triangle $\gamma(x,y) \cup \gamma(y,w) \cup \gamma(w,x)$ is $\delta$-thin by assumption, there exists a vertex $p' \in \gamma(y,w)$ such that $d_{\infty}(p,p') \leq \delta$. Similarly, noticing that
\begin{center}
$d_{\infty}(p',\gamma(y,z)) \geq d_{\infty}(p,\gamma(y,z))-d_{\infty}(p,p') \geq \frac{s}{2}-1-\delta> \delta$,
\end{center}
since $V_{\lfloor s/2 \rfloor +1}, \ldots, V_s$ separates $p$ and $\gamma(y,z)$, we deduce that there exists a vertex $p'' \in \gamma(z,w)$ satisfying $d_{\infty}(p',p'') \leq \delta$, because the geodesic triangle $\gamma(w,y) \cup \gamma(y,z) \cup \gamma(z,w)$ is $\delta$-thin. On the other hand, $H_1, \ldots, H_r$ separates $p$ and $\gamma(w,z)$ (by $d_{\infty}$-convexity of half-spaces), hence
\begin{center}
$r \leq d_{\infty}(p, \gamma(z,w)) \leq d_{\infty}(p,p'') \leq d_{\infty}(p,p')+d_{\infty}(p',p'') \leq 2 \delta$.
\end{center}
We have proved that either $\# \mathcal{V} \leq 4\delta +2$ or $\# \mathcal{H} \leq 2 \delta$. A fortiori, $\min(\# \mathcal{V},\# \mathcal{H}) \leq 4 \delta +2$. $\square$

\section{Acylindrical actions}

\noindent
The main result of this section is the following criterion, establishing the equivalence between several acylindrical properties of an action. Notice however that the acylindrical property we consider is weaker from acylindrical actions as defined in \cite{arXiv:1304.1246} (see Definition \ref{acylindrical}), because of a lack of uniformity on the constants.

\begin{thm}\label{acyl}
Let $G$ be a group acting on a complete CAT(0) cube complex $X$. Suppose $(X,d_{\infty})$ hyperbolic. The following statements are equivalent:
\begin{itemize}
	\item[(i)] for every $d \geq 0$, there exists $R \geq 0$ such that, for every vertices $x,y \in X$,
\begin{center}
$d_{\infty}(x,y) \geq R \Rightarrow \# \{ g \in G \mid d_{\infty}(x,gx),d_{\infty}(y,gy) \leq d\} <+ \infty$;
\end{center}
	\item[(ii)] there exists $R \geq 0$ such that, for every vertices $x,y \in X$,
\begin{center}
$d_{\infty}(x,y) \geq R \Rightarrow \# \{ g \in G \mid gx=x, gy=y\} <+ \infty$;
\end{center}
	\item[(iii)] there exists $R \geq 0$ such that, for any hyperplanes $J_1,J_2$ separated by at least $R$ pairwise disjoint hyperplanes, $\mathrm{stab}(J_1) \cap \mathrm{stab}(J_2)$ is finite.
\end{itemize}
\end{thm}

\noindent
The following lemma will be needed to prove this theorem:

\begin{lemma}\label{fixedpoint}
Let $G$ be a group acting on a complete CAT(0) cube complex $X$. If $G \curvearrowright (X,d_{\infty})$ has a bounded orbit, then $G$ stabilizes a cube.
\end{lemma}

\noindent
It is a consequence of \cite[Theorem A]{Helly}, stating that a bounded Helly graph without infinite clique contains a clique invariant by any automorphism, since it is proved in \cite[Proposition 2.6]{depth} that the cone-off associated to $(X,d_{\infty})$ is a Helly graph. 

\medskip \noindent
\textbf{Proof of Theorem \ref{acyl}.} The implication $(i) \Rightarrow (ii)$ is clear.

\medskip \noindent
Now, we want to prove $(ii) \Rightarrow (iii)$. According to Theorem \ref{dinfty}, there exists a constant $C \geq 1$ such that any grid of hyperplanes in $X$ is $C$-thin; let $R$ denote the constant given by $(ii)$. Let $J_1,J_2$ be two hyperplanes separated by at least $\max(C+1,R)$ pairwise disjoint hyperplanes, say $V_1, \ldots, V_{r}$. Suppose that $K=\mathrm{stab}(J_1) \cap \mathrm{stab}(J_2)$ is infinite and that the action $K \curvearrowright J_1$ has unbounded orbits (with respect to $d_{\infty}$). Thus, if $x \in N(J_1)$ and $y \in N(J_2)$ are two vertices minimizing the combinatorial distance between $N(J_1)$ and $N(J_2)$, there exists $k \in K$ such that $d_{\infty}(x,kx)>C$, ie., there exist $C+1$ pairwise disjoint hyperplanes $H_1, \ldots, H_{C+1}$ separating $x$ and $kx$. A fortiori, the vertices $kx$ and $ky$ also minimize the combinatorial distance between $N(J_1)$ and $N(J_2)$, so that, according to Proposition \ref{hyperplanseparantcor}, the hyperplanes separating $x$ and $y$, and $kx$ and $ky$, are the same: namely, the hyperplanes separating $J_1$ and $J_2$. As a consequence, any hyperplane separating $x$ and $kx$ necessarily separates $y$ and $ky$, because otherwise it would separate $x$ and $y$ or $kx$ and $ky$, but this is impossible; in particular, the $H_i$'s separate $y$ and $ky$. Therefore, the $H_i$'s and the $V_j$'s define a grid of hyperplanes which is not $C$-thin, contradicting the definition of $C$. Consequently, the action $K \curvearrowright J_1$ has a bounded orbit (with respect to $d_{\infty}$); of course, the same statement holds for $K \curvearrowright J_2$. Now, Lemma \ref{fixedpoint} implies that $K$ contains a finite-index subgroup $K_0$ fixing pointwise two cubes $Q_1 \subset N(J_1)$ and $Q_2 \subset N(J_2)$: this contradicts $(ii)$ since $K_0$ is infinite. 

\medskip \noindent
Finally, we want to prove $(iii) \Rightarrow (i)$. Let $R \geq 1$ denote the constant given by $(iii)$. Let $d \geq 1$ and let $x,y \in X$ be two vertices satisfying $d_{\infty}(x,y) \geq R+2d+2$. Our aim is to prove that the set
\begin{center}
$F= \{ g \in G \mid d_{\infty}(x,gx),d_{\infty}(y,gy) \leq d \}$ 
\end{center}
is finite. Let $\mathcal{W}$ be a collection of $R+2d+2$ pairwise disjoint hyperplanes separating $x$ and $y$. We claim that, for every $g \in F$, the images by $g$ of all but at most $2d$ hyperplanes of $\mathcal{W}$ separate $x$ and $y$. Indeed, obviously any hyperplane of $g \mathcal{W}$ must separate $gx$ and $gy$, so that if such a hyperplane does not separate $x$ and $y$, necessarily it has to separate either $x$ and $gx$ or $y$ and $gy$. Therefore, the hypothesis $d_{\infty}(x,gx),d_{\infty}(y,gy) \leq d$ proves our claim.

\medskip \noindent
Thus, if we let $\mathcal{L}$ denote the set of functions from a subset of $\mathcal{W}$ whose cardinality is at least $R+2$ to the set of hyperplanes separating $x$ and $y$, then any element of $F$ naturally induces an element of $\mathcal{L}$. Therefore, because $\mathcal{L}$ is finite, if $F$ is infinite there exist infinitely many elements $g_0,g_1,g_2, \ldots \in F$ inducing the same element of $\mathcal{L}$; in particular, we get infinitely many elements $g_0^{-1}g_1,g_0^{-1}g_2, \ldots \in G$ stabilizing each hyperplane of a subfamily $\mathcal{W}_0 \subset \mathcal{W}$ of cardinality at least $R+2$. Finally, $\mathcal{W}_0$ contains two hyperplanes $J_1,J_2$ separated by at least $R$ pairwise disjoint hyperplanes, such that $\mathrm{stab}(J_1) \cap \mathrm{stab}(J_2)$ is infinite, contradicting $(iii)$. Consequently, $(iii)$ implies that $F$ is necessarily finite. $\square$

\section{Acylindrical hyperbolicity of infinitely-presented groups}\label{small}

\subsection{Small cancellation polygonal complexes}

\noindent
In the following, a polygonal complex $X$ will be a two-dimensional CW-complex whose cells, which are polygons, embed. In this context, a \emph{disc diagram} is a continuous combinatorial map $D \to X$, where $D$ is a finite contractible polygonal complex with a fixed topological embedding into $\mathbb{S}^2$; notice that $D$ may be \emph{non-degenerate}, ie., homeomorphic to a disc, or may be \emph{degenerate}. In particular, the complement of $D$ in $\mathbb{S}^2$ is a $2$-cell, whose attaching map will be referred to as the \emph{boundary path} $\partial D \to X$ of $D \to X$; it is a combinatorial path. 

A disc diagram $D \to X$ is \emph{reduced} if, for any polygons $P,Q \subset D$ meeting along an edge $e$, whenever we think of the boundaries of $P,Q$ as two combinatorial paths $\partial P, \partial Q \to X$ with $e$ as the first edge, then these boundaries are different.

Given a combinatorial closed path $P \to X$, we say that a disc diagram $D \to X$ is \emph{bounded} by $P \to X$ if there exists an isomorphism $P \to \partial D$ such that the following diagram is commutative:
\begin{displaymath}
\xymatrix{ \partial D \ar[rr] & & X \\ P \ar[u] \ar[urr] & & }
\end{displaymath}
According to a classical argument due to van Kampen and Lyndon, there exists a reduced disc diagram bounded by a given combinatorial closed path if and only if this path is null-homotopic.

\begin{definition}
Let $X$ be a polygonal complex. A \emph{piece} in $X$ is a combinatorial path $p$ included into the intersection of two distinct polygons; let $|p|$ denote its length. We say that $X$ satisfies 
\begin{itemize}
	\item the \emph{condition $C'(\lambda)$} if any piece $p$ between two polygons $P,Q$ satisfies $|p| < \lambda |P|$;
	\item the \emph{condition $C(n$)} if the boundary of a polygon cannot be covered by $k$ pieces with $k < n$;
	\item the \emph{condition $T(n)$} if the length of a cycle in the link of a vertex of $X$ is either two or at least $n$.
\end{itemize}
Notice that the condition $C'(1/n)$ implies the condition $C(n)$.
\end{definition}

\noindent
The basic idea is that a reduced disc diagram $D \to X$ in a small cancellation polygonal complex $X$ necessarily contains a polygon with a large part of its boundary in $\partial D$. This is the meaning of the next result.

\begin{definition}
Let $X$ be a polygonal complex and $D \to X$ a disc diagram. A \emph{spur} in $D$ is a vertex of degree one. An \emph{$i$-shell} is a polygon $P$ of $D$ whose boundary can be written as the concatenation of an \emph{outer path} $\partial_{\text{out}} P$ and an \emph{inner path} $\partial_{\text{inn}} P$ such that $\partial_{\text{out}} P = \partial P \cap \partial D$, $\partial_{\text{inn}} P \cap \partial D$ does not contain edges and $\partial_{\text{inn}} P$ is covered by at most $i$ pieces. 
\end{definition}

\noindent
It is worth noticing that a non-degenerate disc diagram containing at least one polygon cannot contain a spur.

\begin{thm}\label{smallcancellationtheorem}\emph{\cite[Theorem 9.4]{McCammondWise}}
Let $X$ be a $C(4)-T(4)$ polygonal complex and $D \to X$ a reduced disc diagram. Then either $D$ contains two spurs and/or $i$-shells with $i \leq 2$, or $D$ is a single vertex or a single 2-cell.
\end{thm}

\noindent
From now on, we will suppose that any polygon of a polygonal complex has an even number of sides. This assumption is not really restrictive since any polygonal complex satisfies it up to a subdivision of its edges; nevertheless, it is worth noticing that this process does not disturb the small cancellation condition. We want to sketch how Wise associates a CAT(0) cube complex to some small cancellation polygon complexes in \cite{MR2053602}.

\begin{definition}
Let $X$ be a polygonal complex. A \emph{hypergraph} $\Lambda$ is an equivalence class of edges with respect to the relation: two edges $e,e'$ are equivalent if there exists a sequence of edges $e_0=e,e_1, \ldots, e_{n-1},e_n=e'$ such that, for every $0 \leq i \leq n-1$, $e_i$ and $e_{i+1}$ are opposite edges in a polygon of $X$. Alternatively, $\Lambda$ may be thought of as a graph included into $X$ whose vertices are the midpoints of the edges which belong to the associated equivalence class and whose edges are segments between two opposite cells in a polygon. The \emph{hypercarrier} of $\Lambda$, denoted by $Y(\Lambda)$, is the union of all the polygons of $X$ intersected by $\Lambda$. 
\end{definition}

\noindent
In \cite{MR2053602}, Wise shows that the hypergraphs and their hypercarriers have good properties. In particular, they embed into the associated polygonal complex, hypercarriers are combinatorially convex, and hypergraphs disconnect the polygonal complex into exactly two connected components. The latter property allows us to define a \emph{halfcarrier} of a hypergraph $\Lambda$ as the union of the hypercarrier $Y(\Lambda)$ of $\Lambda$ with a connected component of the complement of $\Lambda$.

\medskip \noindent
The following description of hypercarriers is especially useful:

\begin{prop}\label{hypercarrier}\emph{\cite[Lemma 3.6]{MR2053602}}
Let $X$ be a $C(4)$ polygonal complex and $\Lambda$ a hypergraph. Then,
\begin{itemize}
	\item the boundary of each polygon of $Y(\Lambda)$ contains a pair of disjoint non-trivial arcs, called \emph{isolated rails}, and the interior of each isolated rail is disjoint from any other polygon of $Y(\Lambda)$;
	\item for each dual edge there is a tree, called a \emph{rung tree}, such that distinct dual edges have disjoint rung trees;
	\item the 1-skeleton of $Y(\Lambda)$ is the union of rung trees and isolated rails such that rung trees corresponding to adjacent vertices of $\Lambda$ are attached together by a pair of isolated rails corresponding to the edge of $\Lambda$ connecting these vertices.
\end{itemize}
\end{prop}

\noindent
Let us mention a consequence of this result which will be useful in the next section.

\begin{cor}\label{corW}
Let $X$ be a $C(4)$ polygonal complex, $R \subset X$ a polygon, $a$ and $b$ two opposite edges in $R$, and $f \nsubseteq R$ an edge adjacent to $R$. Let $Y$ denote the hypercarrier of the hypergraph dual to $a$ and $b$. If $f \subset Y$, then there exists a polygon $S \subset X$ such that either $a$ and $f$, or $b$ and $f$, are included into $\partial S$. 
\end{cor}

\noindent
\textbf{Proof.} If $f \subset Y$, a fortiori there exists a polygon $S$ containing $f$. It follows directly from Proposition \ref{hypercarrier} that the intersection between two polygons in $Y$ is either empty or contains an edge dual to the associated hypergraph. Therefore, $R \cap S \neq \emptyset$ implies that $\partial S$ contains $a$ or $b$. $\square$

\medskip \noindent
As already mentionned, the hypergraphs of a polygonal complex $X$ disconnect the complex into exactly two connected components. In fact, the set of hypergraphs endows $X$ with a structure of a \emph{wallspace}, so that there exists a CAT(0) cube complex $C(X)$ associated to $X$, which is constructed by cubulating this wallspace. See \cite[Section 5]{MR2053602} for more information. In particular, we have a nice description of the maximal cubes in $C(X)$.

\begin{definition}
An \emph{isolated} edge in $X$, ie., an edge which does not belong to any polygon, gives rise to an edge in $C(X)$, which we call an \emph{edge-cube}. A polygon with $2n$ sides in $X$ gives rise to an $n$-cube in $C(X)$, which we call a \emph{cell-cube}.
\end{definition}

\begin{prop}\label{maxcubes}\emph{\cite[Lemma 9.4]{MR2053602}}
Let $X$ be a $C(4)-T(4)$ polygonal complex and $C(X)$ its associated CAT(0) cube complex. Then any cube in $C(X)$ is contained into an edge-cube or into a cell-cube.
\end{prop}

\begin{cor}
The CAT(0) cube complex associated to a $C(4)-T(4)$ polygonal complex is \emph{complete}, ie., it does not contain any infinite increasing sequence of cubes.
\end{cor}

\subsection{A fixed point theorem}

\noindent
In order to link the stabilisers in a small cancellation polygonal complex and in its associated cube complex, we need the following result:

\begin{thm}\label{sc-projection}
Let $G$ be a group acting on a $C'(1/4)-T(4)$ polygonal complex $X$ and let $C(X)$ denote the associated CAT(0) cube complex. There exists a $G$-equivariant projection $p : C(X) \to X$. Moreover, if two points $x,y \in C(X)$ are separated by at least $R+2$ pairwise disjoint hyperplanes, then $p(x)$ and $p(y)$ are separated by at least $R$ pairwise disjoint hypergraphs. 
\end{thm}

\noindent
We begin by proving several preliminary lemmas. From now on, we fix a $C'(1/4)-T(4)$ polygonal complex $X$.

\begin{lemma}
The intersection between two polygons in $X$ is connected.
\end{lemma}

\noindent
\textbf{Proof.} Suppose by contradiction that $X$ contains two polygons $P_1,P_2$ whose intersection is not connected. In particular, there exist two combinatorial paths $\gamma_1 \subset P_1$ and $\gamma_2 \subset P_2$ with the same endpoints, and such that the concatenation $\gamma= \gamma_1 \cup \gamma_2$ defines a simple loop in $X$. Let $D \to X$ be a reduced non-degenerate disc diagram bounded by $\gamma$. If $D$ is a single cell, then its image in $X$ defines a polygon whose boundary is included into $\partial P_1 \cup \partial P_2$: this clearly contradicts the condition $C'(1/4)$. According to Theorem \ref{smallcancellationtheorem}, $D$ must contain a shell; let $R$ denote its image in $X$. The outer path of $R$ is covered by at most two pieces, namely $\partial P_1$ and $\partial P_2$, and its inner path is covered by at most two pieces by definition, so $\partial R$ is covered by at most four pieces: this contradicts the condition $C'(1/4)$. $\square$

\medskip \noindent
The previous lemma will be used frequently in what follows, without any further reference.

\begin{lemma}\label{interpoly}
For any collection of pairwise intersecting polygons $P_1, \ldots, P_n$, the intersection $\bigcap\limits_{i=1}^n P_i$ is non-empty.
\end{lemma}

\noindent
\textbf{Proof.} We argue by induction on $n$. For $n=1,2$, there is nothing to prove. For $n=3$, suppose by contradiction that there exist three pairwise intersecting polygons $P_1,P_2,P_3$ such that $P_1 \cap P_2 \cap P_3= \emptyset$. In particular, there exist three combinatorial paths $\gamma_1 \subset \partial P_1$, $\gamma_2 \subset \partial P_2$, and $\gamma_3 \subset \partial P_3$, such that $\gamma_i \cap \gamma_{i+1}$ is a single vertex $\{v_{i,i+1}\}$ for every $i \in \mathbb{Z}_3$. Thus, the concatenation $\gamma=\gamma_1 \cup \gamma_2 \cup \gamma_3$ defines a simple loop in $X$, and bounds a reduced non-degenerate disc diagram $D \to X$. If $D$ is a single cell, then its image in $X$ is a polygon whose boundary is covered by three pieces, namely $\partial P_1$, $\partial P_2$ and $\partial P_3$: this contradicts the condition $C'(1/4)$. We deduce from Theorem \ref{smallcancellationtheorem} that $D$ contains two shells; let $R_1,R_2$ denote their images in $X$. Because, by definition, the inner path of $\partial R_1$ is covered by at most two pieces, the condition C'(1/4) implies that the outher path of $\partial R_1$ must be covered by at least three pieces: therefore, the outer path of $\partial R_1$ must contain some $\gamma_i$ and intersect the two other subpaths of $\gamma$ along at least one edge. Of course, the same statement holds for $R_2$. But this contradicts the fact that the intersection between the outer paths of $R_1$ and $R_2$ contains no edges. Thus, necessarily $P_1 \cap P_2 \cap P_3 \neq \emptyset$.

\medskip \noindent
Now, let $n \geq 4$ be some integer and suppose that our lemma holds for every $k <n$. By our induction hypothesis, for every $1 \leq k \leq n$, the intersection 
\begin{center}
$I_k= P_1 \cap \cdots \cap P_{k-1} \cap P_{k+1} \cap \cdots \cap P_n$
\end{center}
is non-empty. Suppose by contradiction that $\bigcap\limits_{i=1}^n P_i$ is empty, ie., $I_i \cap I_j = \emptyset$ for every $i \neq j$. For every $k \geq 2$, $I_k$ is a subsegment in $\partial P_1$ whose complement contains $\partial P_k \cap \partial P_1$; notice furthermore that $\partial P_k \cap \partial P_1$ contains $I_j$ for every $j \neq k$. Therefore, $\partial P_1$ is covered by $\partial P_2$, $\partial P_3$, and $\partial P_4$, contradicting the condition C'(1/4). We conclude that $\bigcap\limits_{i=1}^n P_i$ is non-empty. $\square$

\begin{lemma}\label{polyseparated}
Two disjoint polygons are separated by a hypergraph.
\end{lemma}

\noindent
\textbf{Proof.} Let us say that a halfcarrier contains strictly a polygon if this polygon does not belong to the hypercarrier of the associated hypergraph. For any polygon $P$, let $K(P)$ denote the intersection of all the halfcarriers containing strictly $P$ and $\mathrm{star}(P)$ the union of all the closed cells intersecting $P$. To prove our lemma, it is sufficient to show that $K(P)=\mathrm{star}(P)$ for every polygon $P \subset X$. 

\medskip \noindent
The inclusion $\mathrm{star}(P) \subset K(P)$ is clear since $\mathrm{star}(P)$ is included into any halfcarrier containing strictly $P$. Conversely, notice that $K(P)$ is connected, as an intersection of combinatorially convex subcomplexes; therefore, to prove that $K(P) \subset \mathrm{star}(P)$, it is sufficient to show that, for any edge $f \nsubseteq \mathrm{star}(P)$ adjacent to $\mathrm{star}(P)$, we have $f \nsubseteq K(P)$, ie., there exists a hypergraph $\Lambda$ separating $f$ and $P$ such that $f \nsubseteq Y(\Lambda)$. 

\medskip \noindent
\underline{Case 1:} there exists an edge $e$ adjacent to both $P$ and $f$. It follows from Proposition \ref{hypercarrier} that the hypergraph $\Lambda$ dual to $e$ does not intersect $f$ nor $P$: $\Lambda$ separates $f$ and $P$. Let $R$ denote a polygon in $Y(\Lambda)$ containing $e$ and $e' \subset R$ the edge parallel to $e$. According to Corollary \ref{corW}, if $f \subset Y(\Lambda)$ there exists a polygon $S$ containing $e$ and $f$, but this implies $f \subset \mathrm{star}(P)$, or containing $e'$ and $f$, so that $\partial S$ covers at least one half of $\partial R$, contradicting the condition C'(1/4). Therefore, $f \nsubseteq Y(\Lambda)$.

\medskip \noindent
\underline{Case 2:} there exists a polygon $R_0$ adjacent to both $P$ and $f$. Let $R$ be a cell containing $f$ and such that the length of the intersection $|\partial R_0 \cap \partial R|$ is maximal. Let $\gamma_1,\gamma_2$ denote the two subsegments of $\partial R_0$ between $\partial P \cap \partial R_0$ and $\partial R_0 \cap \partial R$. Notice that $f \nsubseteq \mathrm{star}(P)$ implies that the lengths of $\gamma_1$ and $\gamma_2$ are positive.
\begin{center}
\includegraphics[scale=0.6]{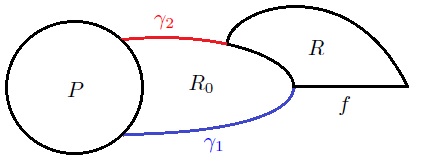}
\end{center}
Suppose by contradiction that no hypergraph intersects both $\gamma_1$ and $\gamma_2$. Let $\overline{\gamma_1}$ denote the subsegment of $R_0$ opposite to $\gamma_1$. Our assumption implies that the intersection $\overline{\gamma_1} \cap \gamma_2$ contains no edges; in particular, $\overline{\gamma_1}$ is included into $\partial R$ or $\partial P$, say $\overline{\gamma_1} \subset \partial R$. Then $\partial R$ covers at least one half of $\partial R_0$, contradicting the condition $C'(1/4)$. Therefore, there exists a hypergraph $\Lambda$ intersecting both $\gamma_1$ and $\gamma_2$. 

\medskip \noindent
If follows from Proposition \ref{hypercarrier} that $\Lambda$ does not intersect $P$ nor $f$, ie., $\Lambda$ separates $P$ and $f$. Let $e_1$ (resp. $e_2$) denote the edge of $\gamma_1 \subset R_0$ (resp. $\gamma_2 \subset R_0$) dual to $\Lambda$. If $f \subset Y(\Lambda)$, Corollary \ref{corW} implies that there exists a polygon $S$ containing either $f$ and $e_1$ or $f$ and $e_2$. If $S$ contains $f$ and $e_1$, then the three polygons $R_0$, $R$ and $S$ define a disc diagram in which the vertex $f \cap R_0$ is interior of degree three, contradicting the condition T(4). If $S$ contains $f$ and $e_2$, then this contradicts the choice of $R$ since 
\begin{center}
$\mathrm{length}(\partial R_0 \cap \partial S) >\mathrm{length}(\partial R_0 \cap \partial R_1)$.
\end{center}
Therefore, $f \nsubseteq Y(\Lambda)$. This concludes the proof. $\square$

\medskip \noindent
\textbf{Proof of Theorem \ref{sc-projection}.} Let $x \in C(X)$. If $x$ belongs to an edge-cube $Q$, there exists an isolated edge $e$ in $X$ with endpoints $y,z \in X^{(0)}$ such that $Q$ is an edge between the principal ultrafilters $v(y)$ and $v(z)$ respectively defined by $y$ and $z$. Now, we define $p(x)$ as the only point of $e$ such that $d(x,v(y))=d(p(x),y)$ and $d(x,v(z))=d(p(x),z)$. 

\medskip \noindent
From now on, suppose that $x$ does not belong to any edge-cube. Thus, according to Proposition \ref{maxcubes}, there is a family of polygons $\{P_i \mid i \in I \}$ naturally associated to the collection of the maximal cubes $\{ Q_i \mid i \in I \}$ containing $x$. It follows from Lemma \ref{polyseparated} that two polygons $P_i$ and $P_j$ intersect if and only if the associated cubes $Q_i$ and $Q_j$ intersect. Therefore, $\{ P_i \mid i \in I \}$ is a collection of pairwise intersecting polygons in $X$. Because the 0-skeleton of a polygon is finite, necessarily there exists a finite subset $J \subset I$ satisfying
\begin{center}
$\bigcap\limits_{i \in I} P_i = \bigcap\limits_{j \in J} P_j$.
\end{center}
It follows from Lemma \ref{interpoly} that this intersection $M(x) \subset X$ is non-empty. If $M(x)$ is a polygon, we define $p(x)$ as the center of this polygon. Otherwise, $M(x)$ is a segment, and we define $p(x)$ as its midpoint. 

\medskip \noindent
The $G$-equivariance of $p : C(X) \to X$ follows by construction. 

\medskip \noindent
Finally, suppose that two points $x,y \in C(X)$ are separated by at least $R+2$ pairwise disjoint hyperplanes. Let $Q(x)$ (resp. $Q(y)$) be a maximal cube containing $x$ (resp. $y$) and let $P(x)$ (resp. $P(y)$) denote the corresponding edge or polygon in $X$. Because the hyperplanes intersecting a cube are pairwise transverse, we deduce that $Q(x)$ and $Q(y)$ are separated by at least $R$ pairwise disjoint hyperplanes; a fortiori, $P(x)$ and $P(y)$ are separated by at least $R$ pairwise disjoint hypergraphs. Since $p(x) \in P(x)$ and $p(y) \in P(y)$, the conclusion follows. $\square$

\medskip \noindent
As a consquence of Theorem \ref{sc-projection}, we are able to deduce the following fixed point theorem:

\begin{cor}\label{FPT}
Let $G$ be a group acting by isometries on a $C'(1/4)-T(4)$ polygonal complex $X$. If $G$ has a bounded orbit, then it has a global fixed point.
\end{cor}

\noindent
\textbf{Proof.} If $G \curvearrowright X$ has a bounded orbit, then the associated action $G \curvearrowright C(X)$ on the CAT(0) cube complex $C(X)$ has also a bounded orbit. Therefore, $C(X)$ contains a global fixed point $x \in C(X)$. Now, if $p : C(X) \to X$ denotes the $G$-equivariant projection given by Theorem \ref{sc-projection}, then $G$ fixes $p(x)$. $\square$

\medskip \noindent
A second fixed point theorem is:

\begin{thm}\label{FPTbis}
Let $G$ be a finitely-generated group acting by isometries on a $C'(1/4)-T(4)$ polygonal complex $X$. If any element of $G$ has a bounded orbit, then $G$ has a global fixed point.
\end{thm}

\noindent
Essentially, this result will be a consequence of the following statement.

\begin{prop}
Let $G$ be a finitely-generated group acting on a CAT(0) cube complex $C$. If $G \curvearrowright (C,d_{\infty})$ contains an unbounded orbit, then there exists a half-space $D$ and an element $g \in G$ such that $gD \subsetneq D$.
\end{prop}

\noindent
A similar statement was proved by Sageev during the proof of \cite[Theorem 5.1]{MR1347406}. We follow his argument below.

\medskip \noindent
\textbf{Proof.} Let $S= \{ s_1, \ldots, s_r \}$ be a finite generating set of $G$ satisfying $S^{-1}=S$ and $v \in C$ a base vertex. By assumption, the orbit $G \cdot v$ is unbounded with respect to $d_{\infty}$, so there exists some $g \in G$ such that $d_{\infty}(v,gv) \geq 3M$ where $M= \sum\limits_{i=1}^r d(v,s_iv)$. Let us write $g$ as a product of elements of $S$: $g=s_{i(1)} \cdots s_{i(m)}$. If, for every $1 \leq i \leq r$, we fix a combinatorial geodesic $\gamma_i$ between $v$ and $s_iv$, then the concatenation
\begin{center}
$\alpha = \gamma_{i(1)} \cup s_{i(1)}\gamma_{i(2)} \cup s_{i(1)}s_{i(2)}\gamma_{i(3)} \cup \cdots \cup s_{i(1)} \cdots s_{i(m-1)}\gamma_{i(m)}$
\end{center}
defines a combinatorial path between $v$ and $gv$.  Now, because $d_{\infty}(v,gv) \geq 3M$, there exists a collection of pairwise disjoint hyperplanes $V_1, \ldots, V_L$ separating $v$ and $gv$, with $L \geq 3M$. Of course, any of these hyperplanes separates $\alpha$, so any $V_i$ is a translate of some $H_j^k$. Because the number of $H_j^k$'s is at most $M$ and $L \geq 3M$, we deduce that there exist three $V_i$'s which belong to the orbit of the same $H_j^k$. 

\medskip \noindent
We have proved that there exist a hyperplane $J$ and two elements $g,h \in G$ such that $J,gJ,hJ$ are three pairwise disjoint hyperplanes separating $v$ and $gv$. For convenience, say that $gJ$ separates $J$ and $hJ$, and let $D$ be the half-space delimited by $J$ which contains $hJ$. If $gD \subset D$ or $hD \subset D$, we have done. Suppose that this not the case, ie., $X \backslash D \subset gD$ and $X \backslash D \subset hD$ or equivalently $D \supset X \backslash gD$ and $D \supset X \backslash hD$. Then 
\begin{center}
$hgD = X \backslash h(X \backslash gD) \subset X \backslash hD \subset D$.  
\end{center}
This concludes the proof. $\square$

\medskip \noindent
\textbf{Proof of Theorem \ref{FPTbis}.} Let $C(X)$ denote the CAT(0) cube complex associated to $X$. Because any element of $G$ has a bounded orbit with respect to the action $G \curvearrowright X$, the same assertion holds with respect to $G \curvearrowright (C(X),d_{\infty})$. Now, the previous proposition implies that the action $G \curvearrowright (C(X),d_{\infty})$ itself has a bounded orbit. Since $C(X)$ is complete, Lemma \ref{fixedpoint} implies that $G$ has a global fixed point $x \in C(X)$. Finally, if $p : C(X) \to X$ denotes the projection given by Theorem \ref{sc-projection}, then $p(x)$ defines a global fixed point in $X$. $\square$

\subsection{Application to small cancellation quotients}

\noindent
Combining the different results established in the previous sections, we are able to find a criterion of acylindrical hyperbolicity for groups acting on a small cancellation polygonal complex. We begin with some preliminaries on acylindrically hyperbolic groups.

\begin{definition}\label{acylindrical}
An action by isometries $G \curvearrowright (S,d)$ is \emph{acylindrical} if, for every $d \geq 0$, there exist two constants $R,N \geq 0$ such that, for every $x,y \in S$,
\begin{center}
$d(x,y) \geq R \Rightarrow \# \{ g \in G \mid d(x,gx),d(y,gy) \leq d \} \leq N$.
\end{center}
A group is \emph{acylindrically hyperbolic} if it admits an action on a (Gromov-)hyperbolic space which is acylindrical and \emph{non-elementary} (ie., with an infinite limit set).
\end{definition}

\noindent
In \cite{arXiv:1304.1246}, Osin gives several equivalent definitions of the acylindrical hyperbolicity of a group. The one we will be interested in is $(AH_3)$ in \cite[Theorem 1.2]{arXiv:1304.1246}. This is because we will use Theorem \ref{acyl} to produce actions on hyperbolic spaces such that any loxodromic isometries will be WPD.

\begin{definition}
Let $G$ be a group acting by isometries on a metric space $(S,d)$. An element $g \in G$ is WPD if, for every $d \geq 0$ and every $x \in S$, there exists $M \geq 1$ such that
\begin{center}
$\# \{ h \in G \mid d(x,hx),d(g^Mx,hg^Mx) \leq d \} < + \infty$.
\end{center}
\end{definition}

\begin{thm}\emph{\cite[Theorem 1.2]{arXiv:1304.1246}}
A group acting by isometries on a hyperbolic space with a WPD isometry is either virtually cyclic or acylindrically hyperbolic.
\end{thm}

\noindent
Our main criterion is the following:

\begin{thm}\label{sc-action}
Let $G$ be a group acting on a $C'(1/4)-T(4)$ polygonal complex $X$. Suppose that there exists a constant $R \geq 0$ such that, for any points $x,y \in X$, the intersection $\mathrm{stab}(x) \cap \mathrm{stab}(y)$ is finite whenever $x$ and $y$ are separated by at least $R$ pairwise disjoint hypergraphs. Either the action $G \curvearrowright X$ is elliptic (ie., any element of $G$ has a fixed a point), or $G$ is virtually cyclic or acylindrically hyperbolic.
\end{thm}

\noindent
The idea is to look at the induced action $G \curvearrowright (C(X),d_{\infty})$ on the CAT(0) cube complex $C(X)$ associated to $X$. Then, we want to show that the small cancellation properties of $X$ implies that $C(X)$ does not contain thick grid of hyperplanes, so that $(C(X),d_{\infty})$ will be hyperbolic, and then use the assumption on the stabilisers in $X$ to deduce an acylindrical property for the action on $C(X)$ by applying Theorem \ref{acyl} thanks to Theorem \ref{sc-projection}. Finally, if we are able to find a loxodromic isometry for the action $G \curvearrowright (C(X),d_{\infty})$, we can deduce that $G$ acts on a hyperbolic space with a WPD isometry. The fist step of this argument is achieved by the following lemma:

\begin{lemma}
Let $X$ be $C'(1/4)-T(4)$ polygonal complex and let $C(X)$ denote its associated CAT(0) cube complex. Then $C(X)$ does not contain $(4,4)$-grids of hyperplanes. 
\end{lemma}

\noindent
\textbf{Proof.} Suppose by contradiction that $C(X)$ contains a $(4,4)$-grid of hyperplanes. In particular, $X$ contains eight hypergraphs $V_1, \ldots, V_4$ and $H_1, \ldots, H_4$ such that $V_i$ and $H_j$ are transverse for every $1 \leq i,j \leq 4$ and $V_i$ (resp. $H_i$) separates $V_{i-1}$ and $V_{i+1}$ (resp. $H_{i-1}$ and $H_{i+1}$) for $i=2,3$. Let $P_{i,j}$ be a polygon in which the hypergraphs $V_i$ and $H_j$ meet. Because $P_{1,1}$ and $P_{1,4}$ belong to the hypercarrier $Y(V_1)$, there exists a sequence of successively adjacent polygons in $Y(V_1)$ between $P_{1,1}$ and $P_{1,4}$; a similar assertion holds for $P_{1,4}$ and $P_{4,4}$, $P_{4,4}$ and $P_{4,1}$, $P_{4,1}$ and $P_{1,1}$. Let $\mathcal{B}$ denote the union of all these polygons, which we call a \emph{band of polygons}. A polygon in $\mathcal{B}$ will be a \emph{corner} if it belongs to $\{ P_{1,1},P_{1,4},P_{4,4},P_{4,1}\}$ and \emph{interior} otherwise. 

\medskip \noindent
Let $\gamma_1 \subset \mathcal{B}$ (resp. $\gamma_2, \gamma_3, \gamma_4 \subset \mathcal{B}$) be a combinatorial path of minimal length between $P_{1,1}$ and $P_{1,4}$ (resp. $P_{1,4}$ and $P_{4,4}$, $P_{4,4}$ and $P_{4,1}$, $P_{4,1}$ and $P_{1,1}$) strictly included into the halfcarrier associated to $V_1$ (resp. $H_4$, $V_4$, $H_1$) which does not contain $V_4$ (resp. $H_1$, $V_1$, $H_4$). Finally, let $c_1$ (resp. $c_2$, $c_3$, $c_4$) be the combinatorial path of minimal length in $P_{1,4}$ (resp. $P_{4,4}$, $P_{4,1}$, $P_{1,1}$) between $\gamma_1$ and $\gamma_2$ (resp. $\gamma_2$ and $\gamma_3$, $\gamma_3$ and $\gamma_4$, $\gamma_4$ and $\gamma_1$). We define the loop $\gamma$ as the concatenation
\begin{center}
$\gamma= \gamma_1 \cup c_1 \cup \gamma_2 \cup c_2 \cup \gamma_3 \cup c_3 \cup \gamma_4 \cup c_4$.
\end{center}
A vertex of $\gamma$ will be called \emph{angular} if it belongs to two polygons of $\mathcal{B}$. 
\begin{center}
\includegraphics[scale=0.6]{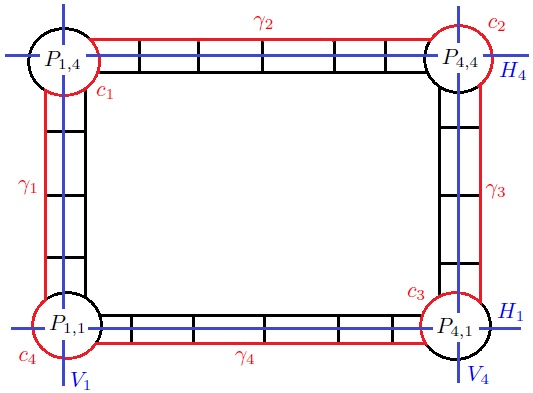}
\end{center}
Notice that the intersection between $c_1$ and $\gamma_1$ or $\gamma_2$ is a single vertex, because we have chosen $\gamma_1$ and $\gamma_2$ of minimal length; furthermore, $c_1$ is disjoint from $c_2$, $c_4$, $\gamma_3$ and $\gamma_4$, because they are separated by $V_2$ or $H_3$. A similar statement holds for $c_2$, $c_3$ and $c_4$. Finally, notice that two $\gamma_i$'s are always separated by a hypergraph so that they are pairwise disjoint. We conclude that $\gamma$ is a simple loop in $X$. Let $D \to X$ be a reduced non-degenerate disc diagram bounded by $\gamma$. 

\medskip \noindent
The reciprocal images in $D$ of the hypergraphs $V_2$, $V_3$, $H_2$ and $H_3$ must separate $D$. In particular, these hypergraphs intersect in different 2-cells in $D$ so that $D$ cannot be a single 2-cell. According to Theorem \ref{smallcancellationtheorem}, $D$ necessarily contains a shell; let $R$ denote its image in $X$. 

\medskip \noindent
We first notice that an interior vertex $v$ of the outer path of $R$ (which is included into $\gamma$) cannot be angular. Otherwise, $v$ belongs to two polygons $P_1,P_2 \subset \mathcal{B}$ and the triple $(P_1,P_2,R)$ defines a disc diagram with $v$ as an interior vertex of degree three, contradicting the condition T(4). Therefore, the outer path of $R$ is included into the boundary of some polygon $P \in \mathcal{B}$. If $P \neq R$, then $\partial R$ is covered by at most three pieces, contradicting the condition C'(1/4). Thus, $R \in \mathcal{B}$.

\medskip \noindent
\underline{Case 1:} $R$ is a corner. Then $R \in \{ P_{1,1},P_{1,4},P_{4,4},P_{4,1}\}$, say $R=P_{1,4}$, and the outer path of $R$ is $c_1$. During the construction of $\gamma$, we have chosen $c_1$ so that $\mathrm{length}(c_1) \leq \frac{1}{2} \mathrm{length}(\partial P_{1,4})$. Consequently, the length of the inner path of $R$ is at least $\frac{1}{2} \mathrm{length}(\partial R)$. By definition, this inner path is covered by at most two pieces, contradicting the condition C'(1/4). 

\medskip \noindent
\underline{Case 2:} $R$ is interior. Then there exist a hypergraph $\Lambda \in \{ V_1,V_4,H_1,H_4\}$ intersecting two edges of the inner path of $R$. Because $\Lambda$ intersects two opposite edges of the polygon $R$, this implies that length of the inner path of $R$ is at least $\frac{1}{2} \mathrm{length}( \partial R)$. By definition, this inner path is covered by at most two pieces, contradicting the condition C'(1/4).

\medskip \noindent
Therefore, we conclude that $C(X)$ does not contain a $(4,4)$-grid of hyperplanes. $\square$

\medskip \noindent
\textbf{Proof of Theorem \ref{sc-action}.} Let $C(X)$ denote the CAT(0) cube complex associated to $X$. According to the previous lemma, $C(X)$ does not contain $(4,4)$-grids of hyperplanes, so we deduce from Theorem \ref{dinfty} that $(C(X),d_{\infty})$ is hyperbolic. Let $x,y \in C(X)$ be two vertices satisfying $d_{\infty}(x,y) \geq R+2$. If $p : C(X) \to X$ denotes the projection provided by Theorem \ref{sc-projection}, then, because $x$ and $y$ are separated by $R+2$ pairwise disjoint hyperplanes, $p(x)$ and $p(y)$ are separated by $R$ pairwise disjoint hypergraphs, so that
\begin{center}
$|\mathrm{stab}(x) \cap \mathrm{stab}(y)| \leq | \mathrm{stab}(p(x)) \cap \mathrm{stab}(p(y)) | < + \infty$.
\end{center}
We deduce from Theorem \ref{acyl} that any loxodromic isometry of $G \curvearrowright (C(X),d_{\infty})$ is WPD. If any element of $G$ has a bounded orbit in $(C(X),d_{\infty})$, then any element has a fixed a point in $C(X)$, and a fortiori in $X$ using Theorem \ref{sc-projection}: thus, the action $G \curvearrowright X$ is elliptic. Otherwise, $G$ has an element $g$ with an unbounded orbit in $(C(X),d_{\infty})$. In particular, $g$ is \emph{combinatorially hyperbolic}, ie., there exists a combinatorial bi-infinite geodesic $\gamma$ on which $g$ acts by translations. By assumption, $\gamma$ has infinite diameter with respect to $d_{\infty}$, so, for every $n \geq 1$, there exist two vertices $x,y \in \gamma$ which are separated by at least $n$ pairwise disjoint hyperplanes. Because $\gamma$ cannot intersect twice a hyperplane, these hyperplanes separate two subrays of $\gamma$; it follows that the Hausdorff distance between two subrays of $\gamma$ is not finite. In particular, $g$ is not parabolic with respect to $d_{\infty}$: it is a loxodromic isometry. 

\medskip \noindent
Therefore, we have proved that $G$ acts on the hyperbolic space $(C(X),d_{\infty})$ with a WPD isometry. We conclude that $G$ is either virtually cyclic or acylindrically hyperbolic. $\square$

\medskip \noindent
The first examples of groups acting on small cancellation polygonal complexes are the small cancellation groups acting on their Cayley complexes. We develop this example below.

\medskip \noindent
Lef $F$ be a free group of rank at least two with a fixed free basis $S$. Every non-trivial element $g$ of $F$ can be written uniquely as a \emph{reduced word} $g=g_1 \cdots g_n$, where each $g_i$ belongs to $S \cup S^{-1}$ and $g_i \neq g_{i+1}^{\pm 1}$. We refer to the integer $n$ as the \emph{length} of $g$, denoted by $|g|$. This reduced word is \emph{cyclically reduced} if $g_1 \neq g_n^{-1}$. If two elements $h,k \in G$, with reduced words $h=h_1 \cdots h_r$ and $k=k_1 \cdots k_s$, satisfy $h_r=k_s^{-1}$, we say that $h_n$ and $k_s$ \emph{cancel} in the product $hk$; otherwise, the product is \emph{reduced}.

Let $\mathcal{R} \subset F$ be a \emph{symmetrised} family, ie., $\mathcal{R}$ is stable by taking cyclic conjugates and inverses. Notice that, up to adding all cyclic conjugates of elements of $\mathcal{R}$ and their inverses, we can always suppose that $\mathcal{R}$ is symmetrised. An element $p \in G$ is a \emph{piece} (with respect to $\mathcal{R}$) if there exist two distinct elements $r_1, r_2 \in \mathcal{R}$ such that the products $r_1=pu_1$ and $r_2=pu_2$ are reduced for some $u_1,u_2 \in G$. 

\begin{definition}
Let $R \subset F$ be a family of cyclically reduced elements and let $\mathcal{R}$ denote its symmetrised. Then $R$ satisfies 
\begin{itemize}
	\item the \emph{condition $C'(\lambda)$} if for every piece $p$ and every element $r \in \mathcal{R}$, such that the product $r=pu$ is reduced, we have $|p| < \lambda |r|$;
	\item the \emph{condition $T(q)$} if, for every $3 \leq h<q$ and every elements $r_0, \ldots, r_{h-1} \in \mathcal{R}$ with $r_i \neq r_{i+1}^{-1}$ for $i \in \mathbb{Z}_h$, at least one of the product $r_1r_2, \ldots, r_{h-2}r_{h-1},r_{h-1}r_1$ is reduced.
\end{itemize}
\end{definition}

\begin{thm}
Infinitely-presented $C'(1/4)-T(4)$ groups are acylindrically hyperbolic.
\end{thm}

\noindent
\textbf{Proof.} Let $G$ be an infinitely-presented $C'(1/4)-T(4)$ group, ie., $G$ admits a presentation $\langle S \mid R \rangle$ where $R$ is an infinite subset of the free group $F(S)$ satisfying the condition $C'(1/4)-T(4)$. Let $X$ denote the Cayley complex associated to this presentation. Notice that $X$ is a $C'(1/4)-T(4)$ polygonal complex, so that $G$ acts freely on a $C'(1/4)-T(4)$ polygonal complex. It follows easily from Theorem \ref{sc-action} that $G$ is either virtually cyclic or acylindrically hyperbolic. On the other hand, as a consequence of Greendlinger Lemma for $C'(1/4)-T(4)$ groups \cite[Theorem V.4.4]{LS}, no relation of $R$ may be deduced from the others, so $G$ is not finitely-presented, and a fortiori $G$ cannot be virtually cyclic. $\square$

\begin{remark}
Finitely-presented $C'(1/4)-T(4)$ are known to be hyperbolic. The situation is similar for $C'(1/6)$ groups: they are hyperbolic in the finitely-presented case, and acylindrically hyperbolic otherwise \cite{C(7)}.
\end{remark}

\begin{ex}
In \cite{TV}, the following group presentation is introduced:
\begin{center}
$G_{I,k}= \langle a,b \mid (a^nb^n)^k=1, \ n \in I \rangle$,
\end{center}
with $k \geq 1$ and $I \subset \mathbb{N} \backslash \{ 0 \}$ infinite. This presentation always satisfies the condition $T(4)$, and the condition $C'(1/4)$ is satisfied precisely when $k \geq 5$. Therefore, $G_{I,k}$ is acylindrically hyperbolic provided that $k \geq 5$. 
\end{ex}

\noindent
The small cancellation theory we have described has been generalized for graph of groups, where we are also able to deduce an action on a small cancellation polygonal complex. Indeed, if $\mathcal{G}$ is a graph of groups with $R \subset \pi_1(\mathcal{G})$ a family satisfying some small cancellation condition, then $\pi_1(\mathcal{G})$ acts on the Bass-Serre tree $T$ associated to $\mathcal{G}$, and, if $\dot{T}$ denotes the usual cone-off of $T$ over the axes of the conjugates of the elements of $R$, then the quotient $\pi_1(\mathcal{G}) / \langle \langle R \rangle \rangle$ acts on the quotient $\dot{T} / \langle \langle R \rangle \rangle$ which is naturally a small cancellation polygonal complex. Below, we describe the situation for free products.

\medskip \noindent
Lef $G=G_1 \ast \cdots \ast G_m$ be a free product. Every non-trivial element $g$ of $G$ can be written uniquely as an \emph{alternating product} $g=g_1 \cdots g_n$, where each $g_i$ is a non-trivial element of a free factor and no two consecutive $g_i,g_{i+1}$ belong to the same free factor. The integer $n$ is called the \emph{free product length} of $g$, denoted by $|g|$. This alternating product is \emph{weakly cyclically reduced} if $|g| \leq 1$ or $g_1 \neq g_n^{-1}$. If two elements $h,k \in G$, with alternating products $h=h_1 \cdots h_r$ and $k=k_1 \cdots k_s$, satisfy $h_r=k_s^{-1}$, we say that $h_n$ and $k_s$ \emph{cancel} in the product $hk$; otherwise, the product is \emph{weakly reduced}.

Let $\mathcal{R} \subset G$ be a \emph{symmetrised} family, ie., $\mathcal{R}$ is stable by taking cyclic conjugates and inverses. Notice that, up to adding all cyclic conjugates of elements of $\mathcal{R}$ and their inverses, we can always suppose that $\mathcal{R}$ is symmetrised. An element $p \in G$ is a \emph{piece} (with respect to $\mathcal{R}$) if there exist two distinct elements $r_1, r_2 \in \mathcal{R}$ such that the products $r_1=pu_1$ and $r_2=pu_2$ are weakly reduced for some $u_1,u_2 \in G$. 

\begin{definition}
Let $R \subset G$ be a family of weakly cyclically reduced elements and let $\mathcal{R}$ denote its symmetrised. Then $R$ satisfies 
\begin{itemize}
	\item the \emph{condition $C'(\lambda)$} if for every piece $p$ and every element $r \in \mathcal{R}$, such that the product $r=pu$ is weakly reduced, we have $|p| < \lambda |r|$;
	\item the \emph{condition $T(4)$} if the two following conditions are satisfied:
	\begin{itemize}
		\item if $r,s,t \in R$ then at least one of the products $rs,st,tr$ is weakly reduced;
		\item if each of $y_1,y_2,y_3$ is a letter occuring in the alternating products of elements $r,s,t \in R$, then $y_1y_2y_3 \neq 1$ in $G$.
	\end{itemize}
\end{itemize}
If $R$ satisfies the condition $C'(\lambda)$, in order to avoid pathological cases, we will use the convention that $|r| \geq \lambda$ for every $r \in R$.
\end{definition}

\begin{thm}\label{sc-freeproduct}
Let $G=G_1 \ast \cdots \ast G_n$ be a finitely-generated free product and $R \subset G$ a family satisfying the condition $C'(1/4)-T(4)$. Then the quotient $Q= G / \langle \langle R \rangle \rangle$ is either virtually cyclic or acylindrically hyperbolic.
\end{thm}

\noindent
\textbf{Proof.} The group $G$ acts on the Bass-Serre tree $T$ associated to its decomposition as a free product. Let $\dot{T}$ denote the usual cone-off of $T$ over the axes of the conjugates of the elements of $R$. Then the quotient $Q=G / \langle \langle R \rangle \rangle$ naturally acts on $X= \dot{T} / \langle \langle R \rangle \rangle$. Naturally, $X$ has the structure of a polygonal complex, and because $R$ satisfies the condition $C'(1/4)-T(4)$, $X$ turns out to be a $C'(1/4)-T(4)$ polygonal complex. As a direct consequence of the well-known Greendlinger Lemma \cite[Theorem 1]{scsurvey}, we get:

\begin{fact}\label{Greendlinger}
If $w \in G \backslash \{ 1 \}$ satisfies $|w| \leq 2$, then $w \notin \langle \langle R \rangle \rangle$.
\end{fact}

\noindent
As a consequence, the free factors of $G$ embed in the quotient $Q$, so that we may identify them with their images into $Q$.

\medskip \noindent
Now, we want to apply Theorem \ref{sc-action}. We first notice that the stabiliser of an edge $e$ in $X$ is the image of $\langle \langle R \rangle \rangle \mathrm{stab}(e') \subset G$ into $Q$, where the edge $e' \subset T$ is a lift of $e$ in $\dot{T}$; because edge-stabilisers in $T$ are trivial, we deduce that edge-stabilisers in $X$ are trivial. Similarly, if $u$ is a vertex of $X$ and $v$ is a lift in $T$, then the stabiliser of $u$ in $Q$ is the image of $\langle \langle R \rangle \rangle \mathrm{stab}(v) \subset G$ into $Q$, ie., it is a conjugate of a free factor. Finally, the stabiliser of a polygon $P$ in $X$ is the image of $\langle \langle R \rangle \rangle \mathrm{stab}(P') \subset G$ into $Q$, where the cone $P' \subset \dot{T}$ is a lift of $P$ in $\dot{T}$. By construction, the basis of the cone $P'$ corresponds to the axis $\gamma$ of a conjugate $f$ of an element $r \in R$. For every $1 \leq i \leq |r|-1$, let $g_i$ be an element of the oriented stabiliser $\mathrm{stab}^+(\gamma)$ such that the translation length of $g_i$ is $i$ modulo $|r|$; if such an element does not exist, set $g_i=1$. Now, if $g$ is any element of $\mathrm{stab}^+(\gamma)$ such that its translation length is $k$ modulo $|r|$, then the translation length of $g_k^{\pm 1}g$ is a multiple of $|r|$, which is also the translation length of $f$; because edge-stabilisers in $T$ are trivial, we conclude that $g_k^{\pm 1}g$ is a power of $f$, so in particular we have $g \in g_k^{\mp 1} \langle f \rangle$. Consequently, the index of $\langle f \rangle$ in $\mathrm{stab}(\gamma)= \mathrm{stab}(P')$ is at most $4 |r|$. A fortiori, since $\langle f \rangle \subset \langle \langle R \rangle \rangle$, the stabiliser of our polygon $P$ in $X$ has cardinality at most $4 |r|$. We have proved:

\begin{fact}\label{stabilisers}
Edge-stabilisers in $X$ are trivial. The stabiliser of a polygon in $X$ associated to an element $r \in R$ has cardinality at most $4|r|$. The vertex-stabilisers in $X$ correspond to the conjugates of the free factors in $Q$. 
\end{fact}

\noindent
Let $x,y \in X$ be two distinct points. We want to prove that the cardinality of $\mathrm{stab}(x) \cap \mathrm{stab}(y)$ is bounded above by $4|r|$ for some $r \in R$. If $x$ or $y$ is not a vertex, then $\mathrm{stab}(x)$ or $\mathrm{stab}(y)$ is included into an edge-stabiliser or a polygon-stabiliser, and the conclusion follows from the previous fact. Therefore, we may suppose that $x$ and $y$ are two vertices. Let $u$ (resp. $v$) be a lift of $x$ (resp. $y$) in $T$. Notice that, because $x \neq y$, $u$ and $v$ are necessarily distinct; in particular, the intersection $\mathrm{stab}(u) \cap \mathrm{stab}(v)$ in $G$ is trivial. Now the intersection $\mathrm{stab}(x) \cap \mathrm{stab}(y)$ in $Q$ is precisely the image of
\begin{center}
$\langle \langle R \rangle \rangle \mathrm{stab}(u) \cap  \langle \langle R \rangle \rangle \mathrm{stab}(v) \subset G$
\end{center}
in the quotient $Q$. Let $g \in G$ be an element of this intersection. So there exist $r_1,r_2 \in \langle \langle R \rangle \rangle$, $g_1 \in \mathrm{stab}(u)$ and $g_2 \in \mathrm{stab}(v)$ such that $r_1g_1=g=r_2g_2$. We deduce that $g_1g_2^{-1}=r_1^{-1}r_2 \in \langle \langle R \rangle \rangle$. It follows from Fact \ref{Greendlinger} that $g_1=g_2 \in \mathrm{stab}(u) \cap \mathrm{stab}(v)=\{1\}$, hence $g \in \langle \langle R \rangle \rangle$. Thus, $\mathrm{stab}(x) \cap \mathrm{stab}(y)$ is trivial. We have proved:

\begin{fact}\label{weakacyl}
For any two distinct points $x,y \in X$, there exists some $r \in R$ such that $\mathrm{stab}(x) \cap \mathrm{stab}(y)$ has cardinality at most $4|r|$.
\end{fact}

\noindent
To conclude, it is sufficient to prove that the action $Q \curvearrowright X$ is not elliptic. Suppose by contradiction that this is the case. According to Theorem \ref{FPTbis}, $Q$ has a global fixed point. So $Q$ contains a finite-index subgroup fixing a vertex $x \in X$. Let $u$ be a lift of $x$ in $T$. Up to a conjugation, we may suppose that $\mathrm{stab}(u)$ is a free factor, say $G_1$. We deduce that $\langle \langle R \rangle \rangle G_1$ is a finite-index subgroup in $G$. On the other hand, our previous fact implies that $\langle \langle R \rangle \rangle G_1 \cap G_2 = \{1 \}$. Thus, the free product $G$ is finite, a contradiction. $\square$

\begin{ex}
Let us introduce, for every $k , p,r,q,s \geq 1$ and $I \subset \mathbb{N} \backslash \{ 0 \}$ infinite, the following group presentation
\begin{center}
$H_{I,k}= \langle a,b,c,d \mid [(ab)^n,(cd)^n]^k=1, a^{p}=b^{q}=c^r=d^s=1, n \in I \rangle$.
\end{center}
The group $H_{I,k}$ is the quotient of the free product 
\begin{center}
$\langle a \mid a^p=1 \rangle \ast \langle b \mid b^q=1 \rangle \ast \langle c \mid c^r=1 \rangle \ast \langle d \mid d^s=1 \rangle$
\end{center}
by the normal closure of the family $R= \{ [(ab)^n,(cd)^n]^k, \ n \in I \}$. The condition $T(4)$ is satisfied if $p,q,r,s \geq 4$, and the condition $C'(1/4)$ is satisfied if $k \geq 5$. Thus, $H_{I,k}$ is acylindrically hyperbolic provided that $p,q,r,s \geq 4$ and $k \geq 5$. 
\end{ex}

\subsection{A note on universal actions}

\noindent
Given an acylindrically hyperbolic group $G$, an element $g \in G$ is a \emph{generalized loxodromic element} if $g$ is loxodromic with respect to some acylindrical action of $G$ on a hyperbolic space; equivalently, $g \in G$ is a generalized loxodromic element if it is a WPD element with respect to an action of $G$ on a hyperbolic space \cite[Theorem 1.4]{arXiv:1304.1246}. An action of $G$ on a hyperbolic space $S$ is a \emph{universal action} if any generalized loxodromic element of $G$ turns out to be WPD; if $G \curvearrowright S$ is moreover acylindrical, this is an \emph{acylindrical universal action}. In \cite[Question 6.7]{arXiv:1304.1246}, Osin asks whether any finitely-generated group admits an acylindrical universal action. A negative answer was given in \cite{Abbott}, by proving that Dunwoody's inaccessible group does not admit an acylindrical universal action. However, the question is still open for finitely-presented groups.

\medskip \noindent
The action we used for small cancellation quotients satisfies the acylindrical property given by Theorem \ref{acyl}, which implies that any loxodromic isometry turns out to be WPD. Thus, we deduce:

\begin{thm}
Any infinitely-presented $C'(1/4)-T(4)$ group $G$ admits a universal action. Furthermore, the generalized loxodromic elements of $G$ are the infinite-order elements.
\end{thm}

\noindent
\textbf{Proof.} We saw that $G$ acts on the CAT(0) cube complex $C(X)$ dual to the Cayley complex associated to the small cancellation presentation; furthermore, $(C(X),d_{\infty})$ is hyperbolic, and the action $G \curvearrowright (C(X),d_{\infty})$ is free and satisfies the acylindrical property given by Theorem \ref{acyl}. Therefore, we deduce that any infinite-order of $G$ induces a loxodromic isometry of $C(X)$, and a fortiori a WPD isometry. $\square$

\medskip \noindent
Unfortunately, we do not know if Theorem \ref{acyl} may be improved to deduce a true acylindricity in general. However such a generalization is possible in finite dimension (where $d_{\infty}$ may be replaced with the combinatorial distance):

\begin{thm}\label{action acylindrique}
Let $G$ be a group acting on an unbounded hyperbolic CAT(0) cube complex $X$. The following are equivalent:
\begin{itemize}
	\item[(i)] $G \curvearrowright X$ is acylindrical,
	\item[(ii)] there exist two constants $L$ and $R$ such that, for any vertices $x,y \in X$ satisfying $d(x,y) \geq L$, the set $\{ g \in G \mid gx=x, gy=y \}$ has cardinality at most $R$,
	\item[(iii)] there exist $N$ and $K$ such that, for any hyperplanes $J_1,J_2$ separated by at least $N$ hyperplanes, the intersection $\mathrm{stab}(J_1) \cap \mathrm{stab}(J_2)$ has cardinality at most $K$.
\end{itemize}
\end{thm}

\noindent
\textbf{Proof.} Without loss of generality, we may suppose that $X$ is $\delta$-hyperbolic for some $\delta \in \mathbb{N}$. According to Theorem \ref{critère d'hyperbolicité}, there exists a constant $C$ such that any grid of hyperplanes in $X$ is $C$-thin.

\medskip \noindent
The implication $(i) \Rightarrow (ii)$ is clear.

\medskip \noindent
Now we prove $(ii) \Rightarrow (iii)$: we suppose that there exist constants $L$ and $R$ such that, for any $x,y \in X$ satisfying $d(x,y) \geq L$, the set $\{ g \in G \mid gx=x, gy=y \}$ has cardinality at most $R$. Without loss of generality, we may assume $L \geq \max(C+1,\dim(X))$. Let $N \geq \mathrm{Ram}(L)$ and $J_1,J_2$ be two hyperplanes separated by at least $N$ hyperplanes. We want to prove that $\left| \mathrm{stab}(J_1) \cap \mathrm{stab}(J_2) \right|$ is bounded above by a constant which depends only on $R$ and $\dim(X)$. 

\medskip \noindent
According to Lemma \ref{Ramsey}, there exist $L$ pairwise disjoint hyperplanes $V_1,\ldots, V_L$ separating $J_1$ and $J_2$. Let $x \in N(J_1)$ and $y \in N(J_2)$ be two vertices minimizing the distance between the two neighborhoods $N(J_1)$ and $N(J_2)$. 

\medskip \noindent
Let $H$ denote the subgroup $\mathrm{stab}(J_1) \cap \mathrm{stab}(J_2)$ and let $g \in H$. Of course, because $d(x,y)=d(gx,gy)$, $gx$ and $gy$ minimize also the distance between $N(J_1)$ and $N(J_2)$, so it follows from Proposition \ref{hyperplanseparantcor} that exactly the same hyperplanes separate $x$ and $y$, and $gx$ and $gy$, namely those which separate $J_1$ and $J_2$. If $J$ is a hyperplane separating $x$ and $gx$, then it cannot separate $x$ and $y$ or $gx$ and $gy$, so it has to separate $y$ and $gy$, and in particular it has to be transverse to the $V_k$'s. Thus, $V_1,\ldots, V_L$ together with the hyperplanes separating $x$ and $gx$ define a $(L,d_{\infty}(x,gx))$-grid of hyperplanes. Because we have supposed $L \geq C+1$, the definition of $C$ implies $d_{\infty}(x,gx) \leq C$, hence $d(x,gx) \leq C \dim(X)$. Similarly, we prove that $d(y,gy) \leq C \dim(X)$. 

\medskip \noindent
We have proved that the orbits of $H$ on $N(J_1)$ and $N(J_2)$ are bounded. Because these neighborhoods are themself CAT(0) cube complexes, this implies that $H$ fixes a point in each one. Furthermore, $H$ stabilizes the maximal cubes which contain these global fixed points, so that $H$ contains a subgroup $H_0$ of index at most $2 \dim(X)!$ which fixes pointwise these two maximal cubes. In particular, $H_0$ has two global fixed vertices $a \in N(J_1)$ and $b \in N(J_2)$. Because $V_1, \ldots, V_{L}$ separate $J_1$ and $J_2$, we deduce that $d(a,b) \geq L$. By our hypotheses, this implies that $H_0$ has cardinality at most $R$. Therefore,
\begin{center}
$| \mathrm{stab}(J_1) \cap \mathrm{stab}(J_2)| =|H| \leq 2 \dim(X)! \cdot |H_0| \leq 2R \dim(X)!$.
\end{center}
This complete the proof.

\medskip \noindent
Now we prove $(iii) \Rightarrow (i)$: we suppose there exist two constants $N$ and $K$ such that, for any hyperplanes $J_1,J_2$ separated by at least $N$ hyperplanes, the intersection $\mathrm{stab}(J_1) \cap \mathrm{stab}(J_2)$ has cardinality at most $K$. Let $\epsilon \in \mathbb{N} \backslash \{ 0 \}$ and let $x,y \in X$ be two vertices satisfying 
\begin{center}
$d(x,y) \geq \mathrm{Ram}(N+2)+2(C \dim(X)+\epsilon+ \mathrm{Ram}(1+C+\epsilon+8\delta))$. 
\end{center}
We want to prove that the cardinality of the set
\begin{center}
$F = \{ g \in G \mid d(x,gx),d(y,gy) \leq \epsilon \}$
\end{center}
is bounded above by a constant which depends only on $\delta$, $C$, $N$, $K$, $\epsilon$ and $\dim(X)$. Thus, the acylindricity of the action will follow. 

\medskip \noindent
Fix a combinatorial geodesic $[x,y]$ and let $z \in [x,y]$ be its midpoint. For convenience, we will suppose in the following that $z$ is a vertex. Let $p,r \in [x,z]$ and $q,s \in [z,y]$ denote the points of $[x,y]$ defined by 
\begin{center}
$d(p,z)=d(q,z)=C \dim(X)+\epsilon+ \lfloor \mathrm{Ram}(N+2)/2 \rfloor$ 
\end{center}
and 
\begin{center}
$d(r,z)=d(s,z)= \mathrm{Ram}(C+\epsilon+8\delta)+C \dim(X)+ \epsilon+ \lfloor \mathrm{Ram}(N+2)/2 \rfloor$. 
\end{center}
Along the geodesic $[x,y]$, we fix our points in the following order: $x$, $r$, $p$, $z$, $q$, $s$ and finally $y$.

\medskip \noindent
Without loss of generality, we may suppose that $C+\epsilon+8\delta \geq \dim X$. Let $g \in F$. We claim that for all but $2(\epsilon+C \dim(X))$ hyperplanes $J$ separating $p$ and $q$, $gJ$ separates $r$ and $s$. If $J$ is a hyperplane separating $p$ and $q$ such that $gJ$ does not separate $r$ and $s$, three cases may happen.

\medskip \noindent
\underline{Case 1:} the hyperplane $gJ$ does not meet $[x,y]$. Because $J$ separates $x$ and $y$, $gJ$ must separate $gx$ and $gy$. Therefore, since $gJ$ does not separate $x$ and $y$, necessarily $gJ$ separates either $y$ and $gy$ or $x$ and $gx$. But $d(x,gx),d(y,gy) \leq \epsilon$, so that there exist at most $2 \epsilon$ such hyperplanes $J$.

\medskip \noindent
\underline{Case 2:} the hyperplane $gJ$ meets $[s,y]$. Let $u \in [p,q]$ be a vertex adjacent to $J$. Then, because in a $\delta$-hyperbolic space the distance is $8\delta$-convex and that $d(x,gx),d(y,gy) \leq \epsilon$, we deduce that $d(u,gu) \leq \epsilon+8\delta$. Let $v \in [s,y]$ be a vertex adjacent to $gJ$. Because $d(q,s) \geq \mathrm{Ram}(C+\epsilon+8\delta)$, Lemma \ref{Ramsey} implies that there exist $k \geq 1+C+\epsilon+8\delta$ pairwise disjoint hyperplanes $V_1,\ldots, V_{k}$ separating $q$ and $s$; say that $V_i$ separates $V_{i-1}$ and $V_{i+1}$ for all $2 \leq i \leq k-1$. Then, each hyperplane $V_i$ separates either $u$ and $gu$ or $v$ and $gu$; in the latter case, $V_i$ is transverse to $gJ$. Because $d(u,gu) \leq \epsilon+8 \delta$, there are at most $\epsilon+8 \delta$ hyperplanes in the first case. On the other hand, because $V_i$ and $V_j$ are disjoint for any $i \neq j$, we deduce that, if $V_i$ is transverse to $gJ$ for some $i$, then $V_j$ is transverse to $gJ$ for all $j \geq i$. Therefore, the hyperplanes $V_{\epsilon+8 \delta},\ldots, V_{k}$ are transverse to $gJ$. Thus, because $k \geq C+1+ \epsilon+8 \delta$, if there exist $M$ hyperplanes $J$ such that $gJ$ meets $[s,y]$, the images by $g$ of these $M$ hyperlanes (which separate $gx$ and $gy$) together with $V_{\epsilon+8 \delta},\ldots, V_{k}$ define a $(\lfloor M/ \dim X \rfloor,C+1)$-grid of hyperplanes. From the definition of $C$, we deduce that $M \leq C \dim(X)$. 

\medskip \noindent
\underline{Case 3:} the hyperplane $gJ$ meets $[r,x]$. This case is symmetric to the previous one. There are at most $C \dim(X)$ such hyperplanes $J$. 

\medskip \noindent
We have just proved that there exist at most $2(\epsilon+C \dim(X))$ hyperplanes $J$ separating $p$ and $q$ such that $gJ$ does not separate $r$ and $s$, as claimed.

\medskip \noindent
Therefore, if $\mathcal{H}(a,b)$ denotes the set of hyperplanes separating two vertices $a$ and $b$, any element $g \in F$ defines a function $(S_g \subset \mathcal{H}(p,q)) \to \mathcal{H}(r,s)$ where $\mathcal{H}(p,q) \backslash S_g$ has cardinality at most $2(\epsilon+C \dim(X))$. The cardinality of the set $\mathcal{L}$ of these functions is bounded above by a constant $\kappa$ which depends only on $d(p,q)$, $d(r,s)$, $C$ and $\epsilon$; from our choices of $d(p,q)$ and $d(r,s)$, in fact it depends only on $\delta$, $C$, $N$, $\epsilon$ and $\dim(X)$.

\medskip \noindent
If $\# F \geq (K+1) \kappa$, then $F$ contains $K+1$ pairwise distinct elements $g_1, \ldots, g_{K+1}$ inducing the same function of $\mathcal{L}$. Therefore, for all $1 \leq i \leq K+1$, $g_1^{-1}g_i$ stabilizes each hyperplane of a family $S \subset \mathcal{H}(p,q)$ such that $\mathcal{H}(p,q) \backslash S$ has cardinality at most $2(\epsilon+C \dim(X))$. On the other hand, $\# \mathcal{H}(p,q) = d(p,q) \geq \mathrm{Ram}(N+2)+2(\epsilon+C \dim(X))$ hence $\# S \geq \mathrm{Ram}(N+2)$. According to Lemma \ref{Ramsey}, the collection $S$ contains $N+2$ pairwise disjoint hyperplanes; because they all separate $p$ and $q$, it makes sense to claim the two extremal hyperplanes $J_1,J_2$ of this subcollection are separated by at least $N$ hyperplanes, hence $|\mathrm{stab}(J_1) \cap \mathrm{stab}(J_2)| \leq K$ by definition of $N$. A fortiori, the intersection $\bigcap\limits_{J \in S} \mathrm{stab}(J)$ has cardinality at most $K$, so we conclude that there exist $i \neq j$ such that $g_1^{-1}g_i =g_1^{-1}g_j$, that is $g_i=g_j$, a contradiction. 

\medskip \noindent
Therefore, $\# F$ is bounded above by the constant $(K+1) \kappa$ which depends only on $\delta$, $C$, $N$, $K$, $\epsilon$ and $\dim(X)$. $\square$

\begin{remark}
The equivalence $(i) \Leftrightarrow (ii)$ of Theorem \ref{action acylindrique} has been proved recently in \cite{articleMartin} for hyperbolic CAT(0) square complexes. In his article, Martin introduces the terminology \emph{weakly acylindrical action} for an action satisfying the condition $(ii)$, so that weakly acylindrical actions and acylindrical actions turn out to be equivalent notions in our context. This statement is especially useful because it is definitely easier to verify whether an action is weakly acylindrical than acylindrical. Therefore, Theorem \ref{action acylindrique} generalizes \cite[Theorem A]{articleMartin} to higher dimensions.
\end{remark}

\noindent
Thus, if the CAT(0) cube complex dual to the small cancellation polygonal complex is finite-dimensional, Theorem \ref{acyl} may be replaced with Theorem \ref{action acylindrique} in the arguments of the previous section in order to produce an acylindrical action. According to Proposition \ref{maxcubes}, we know that our cube complex is finite-dimensional precisely when the lengths of the elements of the small cancellation family we are considering are uniformly bounded. Of course, with respect to the classical small cancellation conditions, this may happen only for finitely-presented groups, which are hyperbolic. On the other hand, for free products we can state:

\begin{thm}
Let $G$ be a free product of finitely-generated groups which are neither acylindrically hyperbolic nor virtually cyclic, and let $R \subset G$ be a family satisfying the condition $C'(1/4)-T(4)$ whose elements have uniformly bounded free product lengths. Then the quotient $G/ \langle \langle R \rangle \rangle$ admits an acylindrical universal action. 
\end{thm}

\noindent
\textbf{Proof.} The group $G$ acts on the Bass-Serre tree $T$ associated to its decomposition as a free product. Let $\dot{T}$ denote the usual cone-off of $T$ over the axes of the conjugates of the elements of $R$. Then the quotient $Q=G / \langle \langle R \rangle \rangle$ naturally acts on $X= \dot{T} / \langle \langle R \rangle \rangle$. Naturally, $X$ has the structure of a polygonal complex, and because $R$ satisfies the condition $C'(1/4)-T(4)$, $X$ turns out to be a $C'(1/4)-T(4)$ polygonal complex. Let $C(X)$ denote the CAT(0) cube complex associated to $X$. We know that $C(X)$ does not contain $(4,4)$-grids of hyperplanes and it is finite-dimensional. Thus, according to Theorem \ref{critère d'hyperbolicité}, $C(X)$ is hyperbolic. Furthermore, combining Fact \ref{weakacyl} and Theorem \ref{action acylindrique}, we deduce that the action $Q \curvearrowright C(X)$ is acylindrical.

\medskip \noindent
Now, because an element of $Q$ inducing an elliptic isometry of $C(X)$ has to fix a point in $C(X)$, and so fixes a point in $X$ as a consequence of Theorem \ref{sc-projection}, and so stabilises a vertex or an edge or a polygon in $X$, we notice thanks to Fact \ref{stabilisers} that an element of $Q$ is elliptic for the action $Q \curvearrowright C(X)$ if and only if it belongs to a conjugate of a free factor or it has finite order. To conclude our proof, it is sufficient to notice that a free factor in $Q$ cannot contain a generalized loxodromic element, but this follows from the assumption that the free factors are neither acylindrically hyperbolic nor virtually cyclic. $\square$

\begin{ex}
For $k \geq 1$ and $I \subset \mathbb{N} \backslash \{ 0 \}$ infinite, let us introduce the group presentation
\begin{center}
$K_{k,I}= \langle a,b,c,d \mid [a,b]=[c,d]=1, (a^nb^nc^nd^n)^k=1, n \in I \rangle$.
\end{center}
The group $K_{k,I}$ is the quotient of the free product of two free abelian groups
\begin{center}
$\langle a,b \mid [a,b]=1 \rangle \ast \langle c,d \mid [c,d]=1 \rangle$
\end{center}
by the normal closure of the family $R= \{ (a^nb^nc^nd^n)^k, n \in I \}$. The condition $C'(1/4)-T(4)$ is satisfied if $k \geq 5$. Therefore, $K_{k,I}$ is acylindrically hyperbolic and admits an acylindrical universal action whenever $k \geq 5$.
\end{ex}

\addcontentsline{toc}{section}{References}

\bibliographystyle{alpha}
\bibliography{CCCC}

\end{document}